\documentclass[3p]{elsarticle}

\usepackage[colorlinks=true]{hyperref}
\usepackage[english]{babel}
\usepackage[utf8]{inputenc}

\usepackage{amsmath}
\usepackage{mathrsfs}
\usepackage{algorithm,algpseudocode}
\usepackage[figbotcap]{subfigure}
\usepackage{enumitem}
\usepackage{MnSymbol}%
\usepackage{wasysym}%
\usepackage{capt-of}

\usepackage{etoolbox}

\makeatletter
\newlength\appendixwidth
\preto\appendix{\addtocontents{toc}{\protect\patchl@section}}
\newcommand{\patchl@section}{%
  \settowidth{\appendixwidth}{\textbf{Appendix }}%
  \addtolength{\appendixwidth}{1.5em}%
  \patchcmd{\l@section}{1.5em}{\appendixwidth}{}{\ddt}%
}
\makeatother

\newcommand\norm[1]{\left\lVert#1\right\rVert}

\journal{Journal of Computational Physics}

\biboptions{sort&compress}

\begin{document}
\hypersetup{
urlcolor=black
}

\begin{frontmatter}
\title{A p-Multigrid Strategy with Anisotropic p-Adaptation Based on Truncation Errors for High-Order Discontinuous Galerkin Methods}

\author[DMAE,CCS]{Andrés M. Rueda-Ramírez \corref{mycorrespondingauthor}}
\cortext[mycorrespondingauthor]{Corresponding authors:}
\ead{am.rueda@upm.es}

\author[DMAE,CCS]{Juan Manzanero}

\author[DMAE,CCS]{Esteban Ferrer}

\author[DMAE,CCS]{Gonzalo Rubio}

\author[DMAE,CCS]{Eusebio Valero}

\address[DMAE]{ETSIAE-UPM (School of Aeronautics), Universidad Politécnica de Madrid, Plaza de Cardenal Cisneros 3, 28040 Madrid, Spain.}

\address[CCS]{Center for Computational Simulation, Universidad Politécnica de Madrid, Campus de Montegancedo, Boadilla del Monte, 28660 Madrid, Spain.}

\begin{abstract}

High-order discontinuous Galerkin methods have become a popular technique in computational fluid dynamics because their accuracy increases spectrally in smooth solutions with the order of the approximation. However, their main drawback is that increasing the order also increases the computational cost. Several techniques have been introduced in the past to reduce this cost. On the one hand, local mesh adaptation strategies based on error estimation have been proposed to reduce the number of degrees of freedom while keeping a similar accuracy. On the other hand, multigrid solvers may accelerate time marching computations for a fixed number of degrees of freedom. 

In this paper, we combine both methods and present a novel anisotropic p-adaptation multigrid algorithm for steady-state problems that uses the multigrid scheme both as a solver and as an anisotropic error estimator. To achieve this, we show that a recently developed anisotropic truncation error estimator \cite[\textit{A. M. Rueda-Ramírez, G. Rubio, E. Ferrer, E. Valero, Truncation Error Estimation in the p-Anisotropic Discontinuous
Galerkin Spectral Element Method, Journal of Scientific Computing}]{RuedaRamirez2018} is perfectly suited to be performed inside the multigrid cycle with negligible extra cost. Furthermore, we introduce a multi-stage p-adaptation procedure which reduces the computational time when very accurate results are required.

The proposed methods are tested for the compressible Navier-Stokes equations, where we investigate two cases. First, the 2D boundary layer flow on a flat plate is studied to assess accuracy and computational cost of the algorithm, where a speed-up of 816 is achieved compared to the traditional explicit method. Second, the 3D flow around a sphere is simulated and used to test the anisotropic properties of the proposed method, where a speed-up of 152 is achieved compared to the explicit method. The proposed multi-stage procedure achieved a speed-up of 2.6 in comparison to the single-stage method in highly accurate simulations.
\end{abstract}

\begin{keyword}
High-order discontinuous Galerkin, anisotropic p-adaptation, error estimation, p-multigrid, compressible flows.
\end{keyword}

\end{frontmatter}

\tableofcontents
\section{Introduction}

Discontinuous Galerkin (DG) methods have gained increasing popularity in the last decades for solving the compressible and incompressible Navier-Stokes equations \cite{Ferrer2017,Manzanero2018b,Wang2013High,Fraysse2016,Browne2014a,Cockburn1998}. 
The lack of a continuity constraint on element interfaces makes DG methods robust for describing advection-dominated problems when an appropriate Riemann solver is selected, and allows them to handle non-conforming meshes with hanging nodes and/or different polynomial orders efficiently \cite{Riviere2008,Kopriva2002,Ferrer2012}.
This is advantageous for accelerating the computations through adaptation strategies that adjust the element size (h) or the polynomial order (p) locally. Multigrid solvers have also been used to accelerate high-order DG time marching computations for a fixed number of degrees of freedom \cite{Fidkowski2005,Luo2006,Luo2006a,Wang2007,Bassi2009,Nastase2006,Shahbazi2009,Wang2009,Mitchell2010,Botti2017}.
The Discontinuous Galerkin Spectral Element Method (DGSEM) \cite{kopriva2009implementing,Black1999} is a high-order nodal variant of the DG technique on hexahedral meshes that is especially suited for mesh adaptation strategies because, in addition to the mentioned properties, it handles p-anisotropic representations efficiently \cite{Kopriva2002,Kompenhans2016,RuedaRamirez2018}.\\

To fully exploit this feature, we can adapt the mesh locally and anisotropically (both in element size and approximation order), so that the solution captures the flow features of interest at a minimum cost. Local adaptation can be performed by subdividing or merging elements (h-adaptation) or by enriching or reducing the polynomial order in certain elements (p-adaptation). To that end, it is paramount to identify the flow regions that require refinement or coarsening. This has been done historically using three different approaches: the \textit{feature-based adaptation}, the \textit{adjoint-based adaptation}, and the \textit{local error-based adaptation}. A comparison of these three approaches was performed by Fraysse et al. \cite{Fraysse2012} in the context of finite volume approximations and by Kompenhans et al. \cite{Kompenhans2016a} and Naddei et al. \cite{Naddei2018} for high-order DG methods. The key ideas behind the adaptation approaches are:\\

\begin{itemize}
\item The \textit{feature-based adaptation} is the classical approach and consists in refining where high velocity, density or pressure gradients are identified. The main disadvantage of these methods is that there is no direct relation between the adaptation criterion and the numerical errors and thus the accuracy is not easily predictable.\\

\item A second and more sophisticated approach is known as \textit{adjoint-based adaptation}. In this approach, a functional target is defined (e.g. drag or lift) and the adjoint problem is solved in order to obtain a spatial distribution of the functional error, which is then used for adapting the mesh. This technique was originally developed for variational formulations \cite{Estep1995,Hartmann2002,Hartmann2006}, and it has been also implemented successfully for Finite Volume schemes \cite{Venditti2002,Pierce2004}. More recently, Wang and Mavriplis \cite{Wang2009} implemented a non-variational formulation for the error estimates and used it to adapt a DG method. The main drawback of this approach is the high computational cost involved in solving the adjoint problem and the storage requirements needed for saving the error estimators.\\

\item A computationally more efficient alternative is the \textit{local error-based adaptation}, which can be based on any local error estimate. On the one hand, estimations of the local discretization error have been used by Mavriplis \cite{Mavriplis1989,Mavriplis1994} to develop hp-adaptation techniques for the spectral element method. Later, Casoni et al. \cite{EvaCasoniyAntonioHuerta2011} extended her approach to adapt the artificial viscosity in shock capturing discontinuous Galerkin discretizations. On the other hand, the $\tau$-estimation method proposed by Brandt \cite{Brandt1984}, which estimates the local truncation error by injecting a fine grid solution into coarser meshes, has been used for adaptation purposes in low-order schemes \cite{berger1987adaptive,Fraysse2012,Fraysse2013,Fraysse2014,Syrakos2012,Syrakos2006}. Rubio et al. \cite{Rubio2013} extended the $\tau$-estimation approach to high-order methods using a continuous Chebyshev collocation method. Later, Rubio et al. applied it to DGSEM discretizations \cite{Rubio2015}, and studied the quasi-\textit{a priori} truncation error estimation, which allows estimating the truncation error without having fully converged fine solutions. Kompenhans et al. \cite{Kompenhans2016} applied the $\tau$-estimation approach to perform p-adaptation using the Euler and Navier-Stokes equations, and showed that a reduction of the truncation error increases the numerical accuracy of all functionals. Furthermore, Kompenhans et al. \cite{Kompenhans2016a} also compared $\tau$-based to featured based adaption, showing better performance for the former. The adaptation strategy consisted in converging a high order representation (reference mesh) to a specified global residual and then performing a single error estimation followed by a corresponding p-adaptation process. More recently, Rueda-Ramírez et al. \cite{RuedaRamirez2018} developed a new method for estimating the truncation error that is cheaper to evaluate than previous implementations, and showed that it produces very accurate extrapolations of the truncation error, which enables using coarser reference meshes.\\

\end{itemize}

The second methodology used in this work are multigrid algorithms. Multigrid methods were first proposed by Brandt \cite{Brandt1977}, who discovered that the classic iterative methods (also referred to as \textit{smoothers}) eliminate the high-frequency components of the error quickly, but fail to eliminate the low-frequency components efficiently. Therefore, he proposed to use coarser meshes to eliminate low-frequency modes. His approach is known as h-multigrid and has been used extensively in low order methods such as traditional Finite Difference and Finite Volume schemes \cite{Hortmann1990,Leister1992,versteeg2007introduction}. Craig and Zienkiewicz \cite{craig1985multigrid}, and Rønquist and Patera \cite{Ronquist1987} were the first authors working on high-order methods that proposed the use of the polynomial order, $p$, to define the levels of a multigrid scheme, the former for p-finite elements and the latter for nodal spectral elements. After these initial works, the use of multilevel methods spread in the high-order community; initially as p-multigrid methods \cite{Fidkowski2005,Luo2006,Luo2006a,Wang2007,Bassi2009} and more recently as hp-multigrid methods \cite{Nastase2006,Shahbazi2009,Wang2009,Mitchell2010,Botti2017}. Most of these implementations use modal hierarchical shape functions \cite{Nastase2006,Wang2007,Fidkowski2005,Shahbazi2009}, and only a small number of publications focus on nodal-based shape functions \cite{Bassi2009,Fidkowski2005}.\\

Two types of multilevel methods can be found in the literature: linear and nonlinear multigrid methods. The former is \textit{de facto} a linear solver and is usually employed to solve the system of equations obtained from an implicit time integration scheme after linearizing with a Newton or Picard iteration. In this case, the smoother is an iterative method for sparse linear systems \cite{saad2003iterative}. The latter, also known as Full Approximation Scheme (FAS), consists in applying the multigrid directly to the set of nonlinear equations. In such a case, the smoother can be either a time-marching scheme (implicit or explicit), or an iterative method applied to the linearized problem. A comparison of linear and nonlinear multigrid methods for DG discretizations can be found in \cite{Nastase2006}. In our work, we make use of the nonlinear multigrid scheme since, as will be shown, it enables the estimation of the truncation error of coarse representations. \\

In this paper, we build on the work on p-adaptation by Kompenhans et al. \cite{Kompenhans2016} and on $\tau$-estimators by Rueda-Ramírez et al. \cite{RuedaRamirez2018}, and combine them with multigrid solution techniques in order to accelerate the convergence of steady-state solutions of the Navier-Stokes equations using the DGSEM. We use the multigrid scheme both as a solver and as an estimator of the truncation error of anisotropic polynomial representations. To do so, we show that the recently developed truncation error estimator by Rueda-Ramírez et al. \cite{RuedaRamirez2018} is well suited to be evaluated inside an anisotropic p-multigrid cycle with a negligible extra cost. The proposed method results in measured speed-ups of up to 816 in a proposed 2D boundary layer case and of 151 in a 3D study of the flow around a sphere, as compared to a traditional explicit solution method. The coupling of multigrid and p-adaptation also enables to propose a multi-stage adaptation process with increasing order representations which reduces the computational cost when very accurate results are required, resulting in speed-ups of 2.6 with respect to the single-stage adaptation process. To the best of our knowledge, this is the first work on DG that couples an anisotropic p-adaptation technique with multigrid.\\

The rest of the paper is organized as follows: in section \ref{sec:DGSEM}, the discontinuous Galerkin spectral element method is briefly explained. Section \ref{sec:AccelTech} details the acceleration methods the current work builds on; namely, multigrid, p-adaptation based on $\tau$-estimations, and their coupling. We finish section \ref{sec:AccelTech} describing how the coupling of multigrid and p-adaptation enables to introduce new features that speed-up the solution procedure. In section \ref{sec:Results}, we study the performance of the proposed p-adaptation algorithm by means of solving 2D and 3D boundary layer test cases. Finally, the main conclusions are gathered in section \ref{sec:Conclusions}.

\section{The Discontinuous Galerkin Spectral Element Method}\label{sec:DGSEM}
We consider the approximation of systems of conservation laws, 

\begin{equation}\label{eq:NScons}
\mathbf{q}_t + \nabla \cdot \mathscr{F} = \mathbf{s},
\end{equation}
where $\mathbf{q}$ is the vector of conserved variables, $\mathscr{F}$ is the flux dyadic tensor which depends on $\mathbf{q}$ and $\nabla \mathbf{q}$, and $\mathbf{s}$ is a source term. As detailed in \ref{sec:NS}, the compressible Navier-Stokes equations can be represented using equation \eqref{eq:NScons}. Multiplying equation \eqref{eq:NScons} by a test function $\mathbf{v}$ and integrating by parts over the domain $\Omega$ yields the weak formulation:
\begin{equation} \label{eq:weak}
\int _{\Omega} \mathbf{q}_t \mathbf{v} \textrm{d} \Omega - \int _{\Omega} \mathscr{F} \cdot \nabla \mathbf{v} \textrm{d} \Omega + \int_{\partial \Omega} \mathscr{F} \cdot \mathbf{n} \mathbf{v} \textrm{d} \sigma = \int _{\Omega} \mathbf{s} \mathbf{v} \textrm{d} \Omega,
\end{equation}
where $\mathbf{n}$ is the normal unit vector on the boundary $\partial \Omega$. Let the domain $\Omega$ be approximated by a tessellation $\mathscr{T} = \lbrace e \rbrace$, a combination of $K$ finite elements $e$ of domain $\Omega^e$ and boundary $\partial \Omega^e$. Moreover, let $\mathbf{q}$, $\mathbf{s}$, $\mathscr{F}$ and $\mathbf{v}$ be approximated by piece-wise polynomial functions (that are continuous in each element) defined in the space of $L^2$ functions:

\begin{equation}
\mathscr{V}^N = \lbrace \mathbf{v}^N \in L^2(\Omega) : \mathbf{v}^N\vert_{\Omega^e} \in \mathscr{P}^N(\Omega^e) \ \ \forall \ \Omega^e \in \mathscr{T} \rbrace,
\end{equation}
where $\mathscr{P}^N(\Omega^e)$ is the space of polynomials of degree at most $N$ defined in the domain of the element $e$. We remark that the functions in $\mathscr{V}^N$ may be discontinuous at element interfaces and that the polynomial order $N$ may be different in each element and direction. Equation \eqref{eq:weak} can then be rewritten for each element as:

\begin{equation} \label{eq:weak2}
\int _{\Omega^e} {\mathbf{q}^e_t}^N {\mathbf{v}^e}^N \textrm{d} \Omega^e - \int_{\Omega^e} {\mathscr{F}^e}^N \cdot \nabla {\mathbf{v}^e}^N \textrm{d} \Omega^e + \int_{\partial \Omega^e} {\mathscr{F}^*} \left( {\mathbf{q}^e}^N, {\mathbf{q}^-}^N, \mathbf{n} \right)  {\mathbf{v}^e}^N \textrm{d} \sigma^e = \int _{\Omega^e} {\mathbf{s}^e}^N {\mathbf{v}^e}^N \textrm{d} \Omega^e,
\end{equation}
where the superindex ``$e$" refers to the functions as evaluated inside the element $e$, i.e. ${\mathbf{q}^e}^N = {\mathbf{q}^N}\rvert_{\Omega^e}$; whereas the superindex ``$-$" refers to the value of the functions on the external side of the interface $\partial \Omega^e$. The numerical flux function, ${\mathscr{F}^*}$, allows to uniquely define the flux at the element interfaces and to weakly prescribe the boundary data as a function of the conserved variable on both sides of the boundary/interface (${\mathbf{q}^{e}}^N$ and ${\mathbf{q}^{-}}^N$) and the normal vector ($\mathbf{n}$). Multiple choices for the numerical flux functions can be found in the literature \cite{toro2013riemann,Manzanero2018a}. In the present work, we use Roe \cite{Roe1981} as the advective Riemann Solver and Bassi-Rebay 1 \cite{Bassi1997} as the diffusive Riemann solver. We remark that the numerical flux must be computed in a specific manner when the representation is non-conforming \cite{Kopriva2002}. \\

Since $\mathbf{q}^N$, $\mathbf{s}^N$, $\mathbf{v}^N$ and $\mathscr{F}^N$ belong to the polynomial space $\mathscr{V}^N$, it is possible to express them inside every element as a linear combination of basis functions $\phi_n \in \mathscr{P}^N(\Omega^e)$:

\begin{eqnarray}\label{eq:PolExp}
\mathbf{q}\rvert_{\Omega^e} \approx {\mathbf{q}^e}^N = \sum_n \mathbf{Q}^e_n \phi^e_n (\mathbf{x}), & & \ \ \ \ 
\mathbf{s}\rvert_{\Omega^e} \approx {\mathbf{s}^e}^N = \sum_n \mathbf{S}^e_n \phi^e_n (\mathbf{x}), \nonumber \\
\mathbf{v}\rvert_{\Omega^e} \approx {\mathbf{v}^e}^N = \sum_n \mathbf{V}^e_n \phi^e_n (\mathbf{x}), & &
\mathscr{F}\rvert_{\Omega^e} \approx {\mathscr{F}^e}^N = \sum_n \pmb{\mathscr{F}}^e_n \phi^e_n (\mathbf{x}).
\end{eqnarray}

Therefore, equation \eqref{eq:weak2} can be expressed in a discrete form as

\begin{equation} \label{eq:DiscretNSElem}
[\mathbf{M}]^e \frac{\partial \mathbf{Q}^e}{\partial t} + \mathbf{F}^e(\mathbf{Q}) = [\mathbf{M}]^e \mathbf{S}^e,
\end{equation}
where $\mathbf{Q}^e=(\mathbf{Q}^e_1, \mathbf{Q}^e_2, \cdots, \mathbf{Q}^e_n, \cdots)^T$ is the local solution that contains the coefficients of the linear combination for the element $e$; $\mathbf{Q}=(\mathbf{Q}^1,\mathbf{Q}^2, \cdots, \mathbf{Q}^K)^T$ is the global solution that contains the information of all elements; $[\mathbf{M}]^e$ is known as the elemental mass matrix, and $\mathbf{F}^e(\cdot)$ is a nonlinear spatial discrete operator on the element level:
\begin{align}
[\mathbf{M}]^e_{i,j} &= \int_{\Omega^e} \phi^e_i \phi^e_j  \textrm{d} \Omega^e, \\
\mathbf{F}^e(\mathbf{Q})_j &= \sum_i \left[ - \int_{\Omega^e} \pmb{\mathscr{F}}_i^e \cdot \phi^e_i \nabla \phi^e_j \textrm{d} \Omega^e \right] + \int_{\partial \Omega^e} {\mathscr{F}^*}^N \left( \mathbf{Q}^e, \mathbf{Q}^{-}, \mathbf{n} \right)  \phi^e_j \textrm{d} \sigma^e.
\end{align}

Note that the operator $\mathbf{F}^e$ is applied on the global solution, since it is the responsible for connecting the elements of the mesh (weakly). Assembling the contributions of all elements into the global system we obtain

\begin{equation} \label{eq:DiscretNS}
[\mathbf{M}] \frac{\partial \mathbf{Q}}{\partial t} + \mathbf{F}(\mathbf{Q}) = [\mathbf{M}] \mathbf{S}.
\end{equation}

In the DGSEM \cite{kopriva2009implementing}, the tesselation is performed with non-overlapping hexahedral elements of order $N=(N_1,N_2,N_3)$ (independent in every direction) and the integrals are evaluated numerically by means of a Gaussian quadrature that is also of order $N=(N_1,N_2,N_3)$. For complex geometries, it is most convenient to perform the numerical integration in a reference element and transform the results to the physical space by means of a high-order mapping of order $M = (M_1,M_2,M_3)$:

\begin{eqnarray} \label{eq:mapping}
\mathbf{x}^e = \mathbf{x}^e \left( \pmb{\xi} \right) \in \mathscr{P}^M, & & \pmb{\xi} = \left( \xi, \eta, \zeta \right) \in \left[ -1, 1 \right]^3.
\end{eqnarray}

The differential operators can be expressed in the reference element in terms of the covariant ($\mathbf{a}_i$) and contravariant ($\mathbf{a}^i$) metric tensors \cite{kopriva2009implementing}:

\begin{equation}
\mathbf{a}_i = \frac{\partial \mathbf{x}^e}{\partial \xi_i}, \ \ \mathbf{a}^i = \nabla \xi_i, \ \ i = 1, 2, 3.
\end{equation}

Using these mappings, the gradient and divergence operators become

\begin{equation}
\nabla q = \frac{1}{J} \sum_{i=1}^d \frac{\partial}{\partial \xi_i} \left( J\mathbf{a}^i q \right), \ \ \nabla \cdot \mathbf{f} = \frac{1}{J} \sum_{i=1}^d \frac{\partial}{\partial \xi_i} \left( J\mathbf{a}^i \cdot \mathbf{f} \right),
\end{equation}
where the Jacobian of the transformation can be expressed in terms of the covariant metric tensor:

\begin{equation}
J = \mathbf{a}_i \cdot \left( \mathbf{a}_j \times \mathbf{a}_k \right), \ \ \left( i,j,k \right) \ \textrm{cyclic}.
\end{equation}

The covariant vectors can be readily obtained from the mapping (equation \eqref{eq:mapping}). For 2D problems, the contravariant vectors can be obtained with the well-known ``cross product form" \cite{Kopriva2006}. However, for fully 3D problems, the contravariant vectors must be obtained using either the ``conservative curl form" or the ``invariant curl form" \cite{Kopriva2006}. Since in this work we deal with 3D curved meshes, the ``invariant curl form" is selected: 

\begin{equation}
Ja_n^i = - \frac{1}{2} \hat x_i \cdot \nabla_{\xi} \times \left[ \mathbf{I}^N \left( X_l \nabla_{\xi} X_m - X_m \nabla_{\xi} X_l \right) \right]  \ \ \ i=1,2,3, \ n=1,2,3, \ \ (n,m,l) \ \text{cyclic},
\end{equation}
where $\mathbf{I}^N$ is an interpolating operator that converts an arbitrary continuous function into a polynomial expansion (as in equation \eqref{eq:PolExp}).\\

Similarly, in the DGSEM the order of the mapping ($M$ in equation \eqref{eq:mapping}) must be $M_i \le N_i$ for 2D, 2D-extruded and 3D p-conforming representations (subparametric or at most isoparametric mapping) to retain free-stream-preservation \cite{Kopriva2006}, whereas it is limited to $M_i \le N_i/2$ for general 3D p-nonconforming representations (David Kopriva, private communication, April 2018). \\

Furthermore, in the DGSEM the polynomial basis functions ($\phi_n$ in equation \eqref{eq:PolExp}) are tensor product reconstructions of Lagrange interpolating polynomials on quadrature points in each of the Cartesian directions of the reference element:

\begin{equation}
\mathbf{q}^N = \sum_n \mathbf{Q}_n \phi_n (\mathbf{x}) = \sum_{i=0}^{N_1} \sum_{j=0}^{N_2} \sum_{k=0}^{N_3} \mathbf{Q}_{i,j,k} l_i (\xi) l_j (\eta) l_k (\zeta).
\end{equation}

Therefore, $\mathbf{Q}_n=\mathbf{Q}_{i,j,k}$ are simply the nodal values of the solution, and $[\mathbf{M}]$ is a diagonal matrix containing the quadrature weights and the mapping terms. In the present work, we make use of the Legendre-Gauss quadrature points \cite{kopriva2009implementing}.\\

A final remark should be made regarding how the time step is chosen. Since in this paper we make use of explicit time integration schemes, the Courant-Friedrich-Levy (CFL) condition dictates a time step limit \cite{karniadakis2013spectral,hesthaven2007nodal}:\\

\begin{equation}
\Delta t = \min(\Delta t^a, \Delta t^{\nu}),
\end{equation}
where the advective time-step restriction is

\begin{equation}
\Delta t^a \le C^a \left( \norm{\pmb{\mathcal{S}}} \frac{N^2}{h} \right)^{-1},
\end{equation}
and the diffusive time-step restriction is 

\begin{equation}
\Delta t^{\nu} \le C^{\nu} \left( \mu \frac{N^4}{h^2} \right)^{-1},
\end{equation}
where $C^a$ and $C^{\nu}$ are constants that depend on the time integration method, $\pmb{\mathcal{S}} = \mathbf{v} + c$ is the characteristic velocity (with $\mathbf{v}$ the flow velocity and $c$ the speed of sound), $\mu$ is the fluid viscosity, and $h$ is the local mesh size. This quantity is evaluated in every time step on the Gauss points of the domain taking into account the possibility of having anisotropic polynomial orders. The most restrictive $\Delta t$ is always chosen.

\section{Acceleration techniques to converge to steady-state} \label{sec:AccelTech}

In this section, we describe the two methods that will be used in the present work to obtain steady-state solutions of the Navier-Stokes equations; namely, nonlinear multigrid schemes and p-adaptation methods based on truncation error estimators.\\

A common way of obtaining a steady-state solution for an unsteady PDE is to start from an arbitrary initial condition and integrate in time until the system converges to a steady solution. The time-stepping scheme can be either explicit or implicit. Explicit high-order Runge-Kutta methods have been traditionally preferred in DGSEM approximations because they do not require solving large complex nonlinear systems that result from implicit implementations. Furthermore, explicit techniques facilitate the parallelization in multi-core systems \cite{Hindenlang2012}. \\

Note that multigrid methods can be (and have been) adapted to unsteady cases in a straightforward manner \cite{Wang2007,Arnone1995,Birken2012}. On the contrary, p-adaptation methods based on $\tau$-estimators have only been applied to steady-state solutions in the context of high-order methods \cite{Kompenhans2016a,Kompenhans2016}. This issue will be addressed in future studies.\\

\subsection{Nonlinear p-multigrid} \label{sec:Multigrid}
Multigrid methods are techniques used to accelerate the convergence to the solution of large linear or nonlinear problems. They constitute a workaround to the fact that standard iterative solvers tend to reduce the high-frequency contents of the error fast, but fail to reduce the low-frequencies efficiently. For this reason, the iterative procedures are commonly referred to as \textit{smoothers} in the multigrid parlance. \\

In h-multigrid, a sequence of progressively coarsening meshes is used where the iterative solver is employed. Every time the mesh is coarsened, some of the smooth components of the error become oscillatory relative to the mesh sampling. Therefore, further coarse-grid smoothing enhances the convergence rate. The p-multigrid scheme relies on the same notion, but low-order polynomial representations are used as coarse levels, which makes it very appropriate for high order methods. p-Multigrid methods typically use a fixed h-mesh.\\ 

For compactness, we shortly describe the Full Approximation Storage (FAS) nonlinear multigrid algorithm. Further details can be found in \cite{Brandt1984,Mavriplis2002,Nastase2006,Wang2007}. Let us consider the steady-state form ($\partial \mathbf{q} / \partial t = 0$) of our nonlinear problem (equation \eqref{eq:DiscretNS}):

\begin{equation}
[\mathbf{M}]^{-1} \mathbf{F}(\mathbf{Q}) = \mathbf{S},
\end{equation}
and define $\mathbf{A}(\mathbf{Q})=[\mathbf{M}]^{-1}\mathbf{F}(\mathbf{Q})$, to obtain

\begin{equation} \label{eq:NonLinFAS}
\mathbf{A}(\mathbf{Q}) = \mathbf{S},
\end{equation}
and use the temporal discretization as the smoothing technique.\\

We select a third order Williamson's low-storage Runge-Kutta scheme (RK3) \cite{williamson1980low} as the smoothing time-marching scheme, so that the developed multigrid schemes can be compared to the purely explicit RK3 (see section \ref{sec:FlatPlateMG}). After some smoothing sweeps in a mesh with polynomial order $P$, the nonlinear residual equation holds

\begin{equation}\label{eq:nonlinRes}
\mathbf{S}^P-\mathbf{A}^P(\mathbf{\tilde Q}^P) =\mathbf{r}^P,
\end{equation}
where $\mathbf{\tilde Q}^P$ is the approximated solution and $\mathbf{r}^P$ is known as the nonlinear residual. Remember that $P = (P_1,P_2,P_3)$ can be different in each element and direction. Using equation \eqref{eq:NonLinFAS}, equation \eqref{eq:nonlinRes} can be rewritten as

\begin{align}
\mathbf{A}^P(\mathbf{Q}^P)-\mathbf{A}^P(\mathbf{\tilde Q}^P) &=\mathbf{r}^P, \label{eq:FASresidual}\\
\mathbf{A}^P(\mathbf{\tilde Q}^P+\pmb{\epsilon}_{it}^P)-\mathbf{A}^P (\mathbf{\tilde Q}^P) &=\mathbf{r}^P,
\end{align} 
where $\pmb{\epsilon}_{it}^P$ is the iteration error on the mesh $P$. The standard two-level p-multigrid FAS scheme consists in transferring equation \eqref{eq:FASresidual} to a lower polynomial representation of order $N = P - \Delta N$ (coarser grid), and using additional smoothing sweeps there. In the lower-order grid, the smoother now targets lower frequencies than the ones removed on the finer grid. Therefore, solving the residual equation on the coarse grid,

\begin{equation}\label{eq:FASCoarseResEq}
\mathbf{A}^N (\mathbf{Q}^N) - \mathbf{A}^N (\mathbf{\tilde Q}_0^N) = \mathbf{r}^N,
\end{equation}
for $\mathbf{Q}^N$, leads to an improved low frequency approximation of the fine grid problem. This holds if $\mathbf{\tilde Q}_0^N$ and $\mathbf{r}^N$ are transferred (interpolated or projected) from the fine grid:
\begin{align}
\mathbf{\tilde Q}_0^N &= \mathbf{I}_P^N \mathbf{\tilde Q}^P \label{eq:FASCoarseSol}\\
\mathbf{r}^{N}          &= \mathbf{I}_P^N \mathbf{r}^P. \label{eq:FASCoarseRes}
\end{align}

Here, $\mathbf{I}_P^N$ is the restriction operator, an $L_2$ projection to the lower polynomial order. Note that no distinction is made between the solution and residual transfer operators since in this work both the solution and the residual are spanned in the same polynomial space. This is an advantage of our implementation since less storage is needed. It is also important to remark that the $L_2$ projection preserves the energy of the transferred quantities. This is an important difference with the transfer operators that are commonly employed in modal DG \cite{Nastase2006,Wang2007,Shahbazi2009,Fidkowski2005}, which do not conserve energy since only the low-order coefficients are transferred for coarse-grid smoothing and the correction is then injected to the lower modes of the high order representation (the transfer matrices are simply identity matrices with rows or columns appended).\\

The coarse-grid nonlinear problem holds
\begin{equation} \label{eq:FAScoarseProblem}
\mathbf{A}^N (\mathbf{Q}^N) = \mathbf{S}^N,
\end{equation}
where $\mathbf{S}^N$ is an artificial source term that can be obtained combining equations \eqref{eq:FASCoarseResEq}, \eqref{eq:FASCoarseSol} and \eqref{eq:FASCoarseRes}:

\begin{equation}\label{eq:FAScoarseSource}
\mathbf{S}^N = \mathbf{A}^N(\mathbf{I}_P^N \mathbf{\tilde Q}^{P}) + \mathbf{I}_P^N \mathbf{r}^P, 
\end{equation}
which, according to equation \eqref{eq:NonLinFAS}, is the same as

\begin{equation}\label{eq:FAScoarseSource2}
\mathbf{S}^N = [\mathbf{M}]^{-1} \mathbf{F}^N(\mathbf{I}_P^N \mathbf{\tilde Q}^{P}) + \mathbf{I}_P^N \mathbf{r}^P.
\end{equation}

After solving equation \eqref{eq:FAScoarseProblem} for $\mathbf{Q}^N$ using a smoothing procedure, we obtain a low frequency approximation of the iteration error:

\begin{equation}\label{eq:FASIterError}
\pmb{\epsilon}^N_{it} = \mathbf{\tilde Q}^N - \mathbf{\tilde Q}^N_0,
\end{equation}
which is then used to update the solution on the fine grid:

\begin{equation}\label{eq:FASCorrectSol}
\mathbf{Q}^P_{i+1}=\mathbf{Q}^P_{i}+ \mathbf{I}_N^P \pmb{\epsilon}^N_{it}.
\end{equation}

The two-level process described above can be generalized to a multilevel FAS V-Cycle and coded efficiently as a recursive procedure as depicted in algorithm \ref{alg:FASVCycle}. Note that the superindex $c$ now denotes the next coarser multigrid level and that the fine superindexes have been dropped for readability. The multigrid cycle has $N_{MG}$ levels, where $level = 1$ is the coarsest (lowest polynomial order) and $level = N_{MG}$ is the finest (highest polynomial order).

\begin{algorithm}[H]
\caption{FAS: V-Cycle}
\label{alg:FASVCycle}
\begin{algorithmic} 
\State \textbf{Recursive Procedure:} FASVCycle( $\mathbf{\tilde Q}$,$\mathbf{r}$,$level$)

\State \textbf{if} $level < N_{MG}$ \textbf{then} $\mathbf{S}=\mathbf{A}(\mathbf{\tilde Q}) + \mathbf{r}$ \Comment{Find coarse-grid source term (eq \eqref{eq:FAScoarseSource})}

\State $\mathbf{\tilde Q}_0 \gets \mathbf{\tilde Q}$ \Comment{Store fine grid solution}
	
\State $\mathbf{\tilde Q} \gets$ \textit{Smooth}($\mathbf{\tilde Q}$,$\beta_1$) \Comment{Pre-smooth $\beta_1$ times (RK3)}

\If{$level > 1$}     
\Comment{If not on the coarsest level, correct the solution using multigrid}
	\State $\mathbf{\tilde Q}^c \gets \mathbf{I}_f^c \mathbf{\tilde Q}$	\Comment{Restrict solution to coarse grid (eq \eqref{eq:FASCoarseSol})}
	\State $\mathbf{r}^c \gets \mathbf{I}_f^c (\mathbf{S}-\mathbf{A}(\mathbf{\tilde Q}))$ \Comment{Restrict residual to coarse grid (eq \eqref{eq:FASCoarseRes})}
	\State CALL \textit{FASVCycle}( $\mathbf{\tilde Q}^c$, $\mathbf{r}^c$,$level-1$) \Comment{Recursive calling}
	\State $\mathbf{\tilde Q}=\mathbf{\tilde Q} + \mathbf{I}_c^f \pmb{\epsilon}^c_{it}$ 		\Comment{Correct solution using coarse-grid approximation (eq \eqref{eq:FASCorrectSol})}
\EndIf
	
\State $\mathbf{\tilde Q} \gets$ \textit{Smooth}($\mathbf{\tilde Q}$,$\beta_2$) \Comment{Post-smooth $\beta_2$ times (RK3)}

\State \textbf{if} $level < N_{MG}$ \textbf{then} $\pmb{\epsilon}_{it} \gets \mathbf{\tilde Q}_0 - \mathbf{\tilde Q}$ \Comment{Compute iteration error (eq \eqref{eq:FASIterError})}

\end{algorithmic}
\end{algorithm}

In this work, we use $\Delta N = P - N = 1$ as the polynomial order reduction in every coarsening. Therefore, the number of multigrid levels corresponds to the maximum polynomial order of the mesh.

\subsubsection{Multigrid cycling strategy} \label{sec:MGcycling}
\begin{figure}[h]
\subfigure[V-cycle.]{\label{fig:Vcycle} \includegraphics[width=0.4\textwidth]{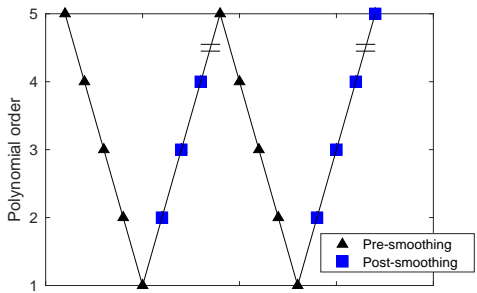}}\qquad
\subfigure[FMG-cycle for getting an appropriate initial condition.]{\label{fig:FMGcycle} \includegraphics[width=0.4\textwidth]{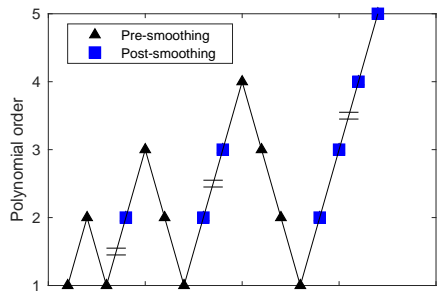}}

\caption{FMG- and V-cycing strategies. Equal signs represent the continuation of the V-cycling process until reaching the desired residual.}\label{fig:MGcycles}

\end{figure}

The typical cycling strategy used for h- and hp-multigrid implementations is to perform repeated V-cycles \cite{Shahbazi2009,Fidkowski2005,Luo2006,Botti2017,Ronquist1987,Nastase2006,Luo2006a,Mitchell2010,Wang2009} (Figure \ref{fig:Vcycle}). Some authors \cite{Shahbazi2009,Mavriplis2002} make use of V or W saw-tooth cycles (without post-smoothing). This technique is well-suited for modal discretizations since the solution correction (equation \eqref{eq:FASIterError}) is injected in the low-order coefficients of the fine-grid representation after coarse-grid smoothing. However, in the nodal discretizations of DGSEM (used in this paper), the coarse-grid smoothed solution can excite high-frequency modes of the fine-grid representation after the interpolation to the fine grid. In consequence, we find that post-smoothing is required.\\

The V-cycling strategy can be very sensitive to the initial condition. To get an appropriate initial condition in the high-order representation, the mainstream alternative is to employ a Full Multigrid (FMG) cycle (see Figure \ref{fig:FMGcycle}) at the beginning of the simulation. Such a cycling strategy can be easily implemented using a recursive algorithm that calls algorithm \ref{alg:FASVCycle} in every ascending level. The number of iterations in every level can be fixed \cite{Shahbazi2009,Nastase2006} or can be tuned using a residual-based approach \cite{Fidkowski2005,Bassi2009}. In this work we use a residual-based approach, where multiple V-cycle repetitions are taken at each level, $p$, until a fixed residual is reached, before raising the approximation level to $p+1$.\\

\subsubsection{Designing the smoothing}\label{sec:TuninSmoothing}
In general, the number of pre-smoothing sweeps ($\beta_1$ in algorithm \ref{alg:FASVCycle}) must be high enough to ensure that the high-frequency modes of the error have been smoothed out, so the $L_2$ restriction does not introduce noise into lower multigrid levels. Likewise, the number of post-smoothing sweeps ($\beta_2$) must be high enough to guarantee that mid-frequency modes of the error do not develop into higher-order representations. These mid-frequency modes of the error can be excited by the $L_2$ prolongation of the solution that was smoothed in a lower multigrid level.\\

A common practice is to set a fixed number of pre- and post-smoothing sweeps \cite{Shahbazi2009,Haupt2013,Bassi2009,Nastase2006,Fidkowski2005}. Nonetheless, when very high polynomial orders and anisotropic non-conforming representations are used, some stages of the simulation can be very sensitive to insufficient smoothing (e.g. at the beginning of the simulation or after an adaptation stage). With that in mind, we propose two residual-based strategies for tuning the number of smoothing sweeps:

\begin{enumerate}
\item \textbf{Pre-smoothing:} After every $\beta_1^0$ sweeps (fixed number), the residual in the next (coarser) representation is checked. If $\norm{\mathbf{r}^{P}}_{\infty} < \eta \norm{\mathbf{r}^{N}}_{\infty}$, the pre-smoothing is stopped; otherwise, $\beta_1^0$ additional sweeps are performed. This strategy is a modification of the residual-based approach that some authors employ in FMG cycles for checking if the coarse level smoothing is enough \cite{Fidkowski2005,Bassi2009}. For the simulations of this paper $\eta \le 1.1$ showed to work fine for meshes with both uniform polynomial orders and also p-anisotropic non-conforming meshes. For this reason, all the simulations that are shown henceforth employ $\eta=1.1$.

\item \textbf{Post-smoothing:} The norm of the residual after the post-smoothing must be at least as low as it was after the pre-smoothing, $\norm{\mathbf{r}^{N}_{post}}_{\infty} \le \norm{\mathbf{r}^{N}_{pre}}_{\infty}$. This condition is checked every $\beta_2^0$ sweeps and the post-smoothing loop is exited when fulfilled. This way, we guarantee that most of the high-frequency errors that could be excited during coarse smoothing are eliminated.
\end{enumerate}

\subsection{p-Adaptation process} \label{sec:AdaptProcess}
Mesh adaptation procedures aim to reduce the number of degrees of freedom of a problem retaining a comparable accuracy. Within those, p-adaptation methods work by increasing the polynomial order of the elements in regions of interest and decreasing it where low order representations are accurate enough. In the present work, we perform anisotropic p-adaptation based on estimations of the truncation error. To do so, we need a methodology for estimating the error of anisotropic polynomial order combinations, which is summarized in next section and detailed in \cite{RuedaRamirez2018}.

\subsubsection{Truncation error estimation} \label{sec:TruncError}
The \textit{non-isolated} truncation error of a discretization of order $N$ ($\tau^N$) is defined as the difference between the discrete partial differential operator ($\mathcal{R}^N$) and the exact partial differential operator ($\mathcal{R}$) applied to the exact solution, $\bar{\mathbf{q}}$: 

\begin{equation}\label{eq:TruncError}
\tau^N = \mathcal{R}^N(\mathbf{I}^N \bar{\mathbf{q}})-\mathcal{R}(\bar{\mathbf{q}}),
\end{equation}
where $\mathbf{I}^N$ is a discretizing operator. For steady-state ($\partial \mathbf{q} / \partial t = 0$), the exact partial differential operator can be derived from equation \eqref{eq:NScons} as

\begin{equation}
\mathcal{R}(\bar{\mathbf{q}}) = \mathbf{s} - \nabla \cdot \mathscr{F} = \bar{\mathbf{q}}_t = 0,
\end{equation}
and the discrete partial differential operator can be derived point-wise from equation \eqref{eq:DiscretNS} as

\begin{equation} \label{eq:Rdiscrete}
\pmb{\mathcal{R}}^N (\pmb{I}^N \bar{\mathbf{q}}) = [\mathbf{M}] \mathbf{S} - \mathbf{F}(\pmb{I}^N \bar{\mathbf{q}}),
\end{equation}
where $\pmb{\mathcal{R}}^N$ contains the sampled values of $\mathcal{R}^N$ in all the nodes of the domain and $\pmb{I}^N$ is a sampling operator. The \textit{non-isolated} truncation error can then be simplified to
 
\begin{equation}\label{eq:TruncErrorSteady}
\pmb{\tau}^N=\pmb{\mathcal{R}}^N(\pmb{I}^N \bar{\mathbf{q}})=[\mathbf{M}] \mathbf{S} - \mathbf{F}(\pmb{I}^N \bar{\mathbf{q}}).
\end{equation}

In addition to the \textit{non-isolated} truncation error, Rubio et al. \cite{Rubio2015} defined the \textit{isolated} truncation error as

\begin{equation}\label{eq:IsolTruncErrorSteady}
\hat{\pmb{\tau}}^N=\hat{\pmb{\mathcal{R}}}^N(\pmb{I}^N\bar{\mathbf{q}})=[\mathbf{M}] \mathbf{S} - \mathbf{\hat F}(\pmb{I}^N\bar{\mathbf{q}}),
\end{equation}
where $\mathcal{\hat R}^N(\cdot)$ is the \textit{isolated} discrete partial differential operator, which is derived without substituting the flux, $\mathscr{F}$, by the numerical flux ,$\mathscr{F}^*$, in equation \eqref{eq:weak2}, thus eliminating the influence of neighboring elements and boundaries in the truncation error of each element ( see \cite{Rubio2015,RuedaRamirez2018}). Rubio et al. \cite{Rubio2015} showed that the \textit{isolated} truncation error in an element depends on the polynomial order of the element, whereas the \textit{non-isolated} truncation error in the same element depends on the polynomial order of the element \textit{and} its neighbors. In consequence, it was suggested that the \textit{isolated} truncation error may be a better sensor for local p-adaptation than the \textit{non-isolated} truncation error since it is not contaminated by neighbors' errors. Moreover, Rueda-Ramírez et al. \cite{RuedaRamirez2018} showed that the accurate estimation of the \textit{isolated} truncation error imposes fewer conditions, and therefore can be computationally cheaper, than the one of its \textit{non-isolated} counterpart. In this work, we retain the isolated truncation as the driver of the proposed p-adaptation procedure.\\

The aim of this work is the development of a method for solving the Navier-Stokes equations in complex geometries, where the exact solution, $\bar{\mathbf{q}}$, is usually not at hand. Therefore, we utilize the $\tau$-estimation method, which approximates the truncation error using an approximate solution on a high order grid, $P$, instead of the exact one. Furthermore, we are interested in a low cost approximation which suits the multigrid procedure. In consequence, we use the \textit{quasi a-priori} approach without correction \cite{Kompenhans2016}, which makes use of a non-converged solution, $\mathbf{\tilde Q}^P$:

\begin{equation}\label{eq:TruncError!}
\pmb{\tau}_P^N = [\mathbf{M}^N]\mathbf{S}^N-\mathbf{F}^N(\mathbf{I}_P^N \mathbf{\tilde Q}^P).
\end{equation}

Here, the estimate of the \textit{isolated} truncation error can be obtained by simply replacing $\mathbf{F}^N$ by $\hat{\mathbf{F}}^N$ in equation \eqref{eq:TruncError!}. In the rest of this work, the expressions containing the symbol $\tau$ are valid for both the \textit{non-isolated} and the \textit{isolated} truncation error unless the contrary explicitly stated. Kompenhans et al. \cite{Kompenhans2016} showed that $\mathbf{\tilde Q}^P$ must be converged down to a residual $\tau_{max}/10$, in order for equation \eqref{eq:TruncError!} to yield accurate estimations of $\tau^N$ in regions where $\tau^N > \tau_{max}$. In addition, Rueda-Ramírez et al. \cite{RuedaRamirez2018} showed that for p-anisotropic representations, the truncation error of a polynomial order combination, $N=(N_1,N_2,N_3)$, can be obtained as the sum of individual directional components:

\begin{equation} \label{eq:AnisTruncError}
\tau^{N_1N_2N_3} \approx \tau_1^{N_1N_2N_3} + \tau_2^{N_1N_2N_3} + \tau_3^{N_1N_2N_3} \approx \tau_{P_1P_2P_3}^{N_1P_2P_3} + \tau_{P_1P_2P_3}^{P_1N_2P_3} + \tau_{P_1P_2P_3}^{P_1P_2N_3}.
\end{equation}

Each of the directional components, $\tau_i$, has spectral convergence with respect to the polynomial order in the corresponding direction, $N_i$ \cite{RuedaRamirez2018}. This allows obtaining accurate extrapolations of $\tau$ by extrapolating the values of $\tau_i$ with a \textit{linear-log} regression and summing the individual contributions \cite{RuedaRamirez2018}.

\subsection{Coupling anisotropic $\tau$-estimation-based adaptation with multigrid}\label{sec:Coupling}

In this section, we present a new technique for obtaining steady-state solutions based on coupling anisotropic p-adaptation methods and multigrid. As pointed out firstly by Brandt \cite{Brandt1984}, and then recently used by Syrakos et al. \cite{Syrakos2012} in the context of h-refinement techniques, the concept of truncation error arises naturally in FAS multigrid methods. In fact, the second term of our \textit{non-isolated} truncation error estimator (equation \eqref{eq:TruncError!}) is contained in the coarse grid source term of the multigrid scheme (equation \eqref{eq:FAScoarseSource2}). Consequently, computing $\tau_P^N$ inside the multigrid cycle only involves a few additional operations. In the case of the \textit{isolated} truncation error, an additional inexpensive step is required to evaluate the operator $\hat{\mathbf{F}}^N$.\\

Two main differences with previous $\tau$-estimators can be identified. First, instead of interpolating the finest solution directly to every coarser representation, the solution is interpolated level by level in multigrid methods; and second, the smoothing procedure modifies the finest solution before it is transferred to lower orders. In that regard, preliminary tests showed no significant difference between the multigrid $\tau$-estimations and the conventional ones. \\

Since in p-multigrid techniques the coarsening is usually performed in all directions simultaneously, i.e. the polynomial order of every direction is decreased (isotropic multigrid), only certain combinations of low-order polynomial orders are evaluated. This makes it impossible to generate the full anisotropic truncation error map \cite{RuedaRamirez2018} (that is needed for performing anisotropic p-adaptation) with the conventional $\tau$-estimation procedure. For this reason, in section \ref{sec:AnisFAS} we propose a p-anisotropic multigrid procedure and explain how the anisotropic decoupled truncation error estimator by Rueda-Ramírez et al. \cite{RuedaRamirez2018} can be evaluated using such a multigrid scheme for generating the full truncation error map. In section \ref{sec:NewPolorders}, we describe how the new polynomial orders are computed based on the proposed estimations; and finally, we present a multi-stage p-adaptation process in section \ref{sec:MultiStage}. \\

\subsubsection{Anisotropic multigrid} \label{sec:AnisFAS}
The classical approach to implement a p-multigrid method is to perform coarsening in all directions simultaneously. This strategy will be referred to as isotropic multigrid in following sections. In this paper, we propose the use of an anisotropic multigrid method in which the coarsening is done in each direction at a time, in order to estimate the truncation error. For instance, in a 3D p-anisotropic multigrid case, when coarsening in $\xi$, the coarse grid problem is derived from equation \eqref{eq:FAScoarseProblem} as

\begin{equation} 
\mathbf{A}^{N_1P_2P_3} (\mathbf{Q}^{N_1P_2P_3}) = \mathbf{S}^{N_1P_2P_3},
\end{equation}
where the source term is obtained from equation \eqref{eq:FAScoarseSource2}:

\begin{equation} 
\mathbf{S}^{N_1P_2P_3} = [\mathbf{M}]^{-1} \mathbf{F}^{N_1P_2P_3}(\mathbf{I}_{P_1P_2P_3}^{N_1P_2P_3} \mathbf{\tilde Q}^{P_1P_2P_3}) + \mathbf{I}_{P_1P_2P_3}^{N_1P_2P_3} \mathbf{r}^{P_1P_2P_3}.
\end{equation}

If the anisotropic p-multigrid method is used to solve a p-anisotropic representation, the number of multigrid levels can be different in every direction, $N_{MG,i}$.\\

\begin{figure}[h]
\begin{center}
\includegraphics[width=0.6\textwidth]{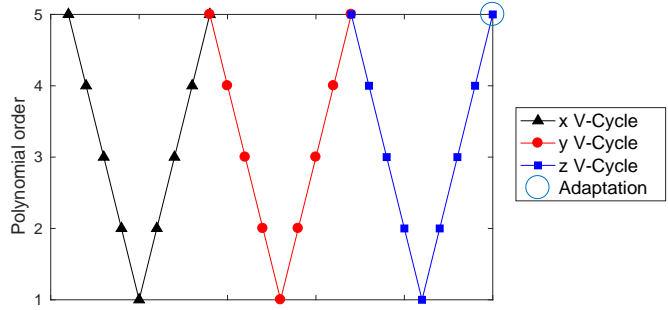}
\caption{Adaptation process: Anisotropic 3V FAS cycle and subsequent adaptation.}\label{fig:AnisFAS}
\end{center}
\end{figure}

Note that this method is perfectly suited to generate the truncation error map using the decoupled truncation error estimator proposed by Rueda-Ramírez et al. \cite{RuedaRamirez2018} (equation \eqref{eq:AnisTruncError}). Figure \ref{fig:AnisFAS} depicts the so-called anisotropic 3V FAS cycle. In every V-cycle, the coarsening is performed in one of the coordinate directions of the reference element and the directional component of the truncation error is estimated. Afterwards, the p-adaptation process that is detailed in section \ref{sec:NewPolorders} takes place. It is noteworthy that the reference coordinate frame of an element inside a general 3D mesh is commonly not aligned with its neighbors'. This can pose a problem for the \textit{non-isolated} truncation error estimation, but not for the \textit{isolated} truncation error that neglects the contribution of the neighboring elements (a thorough analysis can be found in \cite{RuedaRamirez2018}).  \\

Notice that instead of evaluating every possible combination of $N=(N_1,N_2,N_3)$ for $N_{i}<P_i$ (which can be a large number for 3D cases), as in \cite{Rubio2015,Kompenhans2016},  the full truncation error map is constructed from a completely decoupled approach. One important advantage of doing so is that all the storage needed for the decoupled error estimators is already allocated in the anisotropic multigrid routines and only a few inexpensive  additional computations are required. Hence, the multigrid process works indeed as both solver and $\tau$-estimator.\\

\subsubsection{Computing the new polynomial orders} \label{sec:NewPolorders}
Given a truncation error threshold, $\tau_{max}$, that needs to be achieved in a specific case, and a maximum polynomial order allowed, $N_{max}$, the proposed adaptation process can be summarized in six steps:

\begin{enumerate}
\item A high-order representation, $P= (P_1,P_2,P_3)$, is converged down to a residual $\tau_{max}/10$ using the multigrid method described in section \ref{sec:Multigrid}.
\item An anisotropic multigrid procedure (section \ref{sec:AnisFAS}) is used to estimate the decoupled truncation error contribution of every direction. For instance, when coarsening in the direction $(1)$, the contribution is:
\begin{equation}
\pmb{\tau}_{1}^{N_1N_2N_3} \approx \pmb{\tau}_{P_1P_2P_3}^{N_1P_2P_3} = [\mathbf{M}^{N_1P_2P_3}] \mathbf{S}^{N_1P_2P_3}-\mathbf{F}^{N_1P_2P_3}(\mathbf{I}_{P_1P_2P_3}^{N_1P_2P_3}\tilde{\mathbf{Q}}^{P_1P_2P_3}).
\end{equation}

\item The \textit{inner} truncation error map (for $N_i<P_i$) is generated using equation \eqref{eq:AnisTruncError}.

\item If $\tau_{max}$ can be achieved using one of these combinations, it is selected and the simulation continues.
\item If $\tau_{max}$ cannot be achieved in the \textit{inner} truncation error map, an extrapolation procedure based on linear-log regression is conducted in each of the three directions of the decoupled truncation error ($\tau_i$), and then the full truncation error map is generated for $P_i \le N_i \le N_{max,i}$. 

\item If $\tau_{max}$ can be achieved using one of these combinations, it is selected. If not, $N_1 = N_2 = N_3 = {N}_{max}$ is selected.
\end{enumerate}

In steps 4 and 6, there can be multiple combinations ($N_1$,$N_2$,$N_3$) that achieve $\tau < \tau_{max}$. In that case, the combination with the lowest number of degrees of freedom is selected. Notice that the two main differences with the method of Kompenhans et al. \cite{Kompenhans2016} are: ($i$) the way in which the truncation error is estimated for $N_i<P_i$ (step 3) and later for $N_i \ge P_i$ (step 5); and ($ii$) that if the truncation error is not achieved, the element is fully enriched in all directions, instead of in only one.

\subsubsection{Uniform coarsening versus high-order coarsening} \label{sec:HOcoarsening}
We propose two ways of obtaining the polynomial orders of the coarser representations. 
Let us now define the diadic tensor $\mathcal{N}$, which contains the polynomial orders of all the elements in a mesh:

\begin{equation}
\mathcal{N} = (N^1, N^2, \cdots, N^{e}, \cdots, N^K),
\end{equation}
where $K$ is the number of elements of the mesh, $e$ is the element index, and $(N^{e}=N^{e}_1,N^{e}_2,N^{e}_3)$.\\

After the p-adaptation procedure is done, the mesh consists of elements with non-uniform anisotropic polynomial orders. Taking into account that $\Delta N = 1$ (section \ref{sec:Multigrid}), the number of multigrid levels is $N_{MG} = \max(\mathcal{N}) - N_{\textit{coarse}} + 1$ for the isotropic multigrid and $N_{MG,i} = \max(\mathcal{N}_i) - N_{\textit{coarse}} + 1$ for the anisotropic multigrid (the latter is a function of the maximum polynomial order per direction). Let us define two ways of performing the coarsening inside a multigrid cycle: 

\begin{itemize}
\item \textit{Uniform coarsening}: The coarsening is performed in all elements simultaneously:
\begin{equation} \label{eq:coarsening_ini}
\left( N^{e}_i \right)_{level} = \left( N^{e}_i \right)_{level+1} - \Delta N,
\end{equation} 
except in the elements where the minimum polynomial order allowed has been reached:
\begin{equation}
if \ \left( \left( N^{e}_i \right)_{level}<N_{coarse} \right) \ \mathbf{then} \ \left( N^{e}_i \right)_{level}=N_{coarse}.
\end{equation}

\item \textit{High-order coarsening}: Since the maximum polynomial order in every multigrid level can be known beforehand:
\begin{equation} \label{eq:highorderCoarse}
\left( N_i \right)_{level} ^{max} = \max(\mathcal{N}_i) - \Delta N (N_{MG,i} - level),
\end{equation}
we can coarsen only the elements that do not fulfill this condition:
\begin{equation} \label{eq:coarsening_fin}
if \ ( \left( N_i^{e} \right)_{level+1} > \left( N_i \right)_{level}^{max}) \ \mathbf{then} \ \left( N_i^{e} \right)_{level}  = \left( N_i \right)_{level}^{max}.
\end{equation}
In this way only the high-order elements are coarsened. In this paper, we use $N_{\textit{coarse}}= \Delta N = 1$. Therefore, equation \eqref{eq:highorderCoarse} reduces to
\begin{equation}
\left( N_i \right)_{level}^{max} = level.
\end{equation}
\end{itemize}

Notice that equations \eqref{eq:coarsening_ini} to \eqref{eq:coarsening_fin} are valid for isotropic and anisotropic multigrid procedures. In the former, $N_{MG,i}$ must be simply replaced by $N_{MG}$, $\max(\mathcal{N}_i)$ by $\max(\mathcal{N})$, and the operations are performed in all directions. In the latter, the operations are only performed in the direction in which the coarsening is done. Furthermore, both coarsening methods are equivalent for isotropic polynomial representations.\\

One could argue that the \textit{uniform coarsening} involves less computational cost than the \textit{high-order coarsening} since coarse representations have fewer degrees of freedom. Nevertheless, the latter has two main advantages:

\begin{itemize}
\item Several preliminary tests showed that \textit{uniform coarsening} could be unstable for highly anisotropic meshes in 2D and 3D. 

\item In 3D meshes (that are not 2D extrusions), p-nonconforming representations require the mapping order to be $M \le N/2$ (see section \ref{sec:DGSEM}). This means that the minimum polynomial order of the mesh must be $\min(\mathcal{N}) \ge 2$. If the \textit{uniform coarsening} is used, the mapping restriction forces the coarsest multigrid level to have a polynomial order $N_{\textit{coarse}} \ge 2$. However, if the \textit{high-order coarsening} is used, the coarsest polynomial order can be as low as $N_{\textit{coarse}} \ge 1$ since the two coarsest levels are always p-conforming. The additional coarse multigrid level helps to eliminate the low frequency components of the error.
\end{itemize}

For these reasons, in this paper we use only \textit{high-order coarsening}.

\subsubsection{Multi-stage adaptation process} \label{sec:MultiStage}
\begin{figure}[h]
\begin{center}
\includegraphics[width=0.6\textwidth]{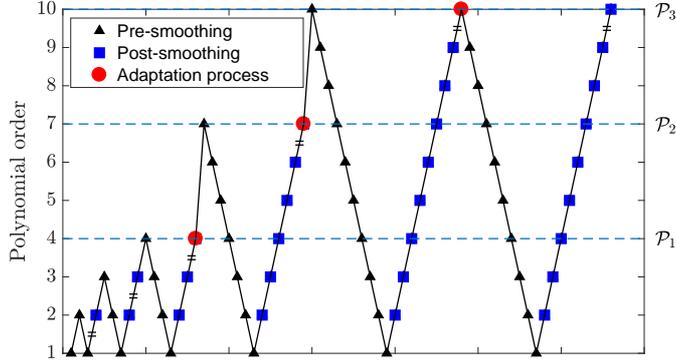}
\caption{Proposed FMG cycle with multiple adaptation stages. Equal signs represent the continuation of the V-cycling process until reaching the desired residual (see section \ref{sec:MGcycling}).}\label{fig:FMGAdapt}
\end{center}
\end{figure}

The proposed \textit{multi-stage} adaptation strategy takes advantage of an FMG-cycle and performs multiple adaptation processes at different polynomial orders ($\mathcal{P}_i$), as depicted in Figure \ref{fig:FMGAdapt}. In an adaptation stage at level $\mathcal{P}_i$ (red circular markers), a $\tau$-estimation procedure is performed using a $3V$ anisotropic multigrid cycle (Figure \ref{fig:AnisFAS}), and subsequently, the polynomial orders are adjusted accordingly, but never to a polynomial order that is higher than $\mathcal{P}_{i+1}$.  In such a case, $\mathcal{P}_{i+1}$ is selected. This differs from the previous adaptation strategies based on $\tau$-estimation in that, traditionally, the whole domain had to be solved in a considerably high-order mesh before performing the single-stage adaptation process. The main advantage of using this methodology is that the zones of the domain that only require a low order representation are identified early in the simulation and are not enriched. This reduces the overall computational costs.\\

After a truncation error estimation at the level $\mathcal{P}_i$, the algorithm checks if the maximum \textit{required} polynomial order is lower or equal to the polynomial order of the next stage, $\mathcal{P}_{i+1}$. In such a case, the performed adaptation step is marked as the last one and the simulation continues without any further adaptation processes.

\section{Numerical results} \label{sec:Results}
In this section, we test the accuracy and performance of the proposed p-adaptation algorithms in 2D (section \ref{sec:FlatPlate}) and 3D (\ref{sec:Sphere}) cases. For the reasons exposed in section \ref{sec:TruncError}, we will use the \textit{isolated} truncation error for the p-adaptation algorithms. All the results presented in this section were obtained using an 8-core 2.6 GHz Intel Xeon E5-2670 and 32 GB RAM, and shared memory parallelization (OpenMP) for computing the spatial terms, as explained in \cite{Hindenlang2012}. Note that our parallel implementation has not been optimized for a specific cache architecture and further speed-ups may be achieved.\\

\subsection{2D Flow over a flat plate} \label{sec:FlatPlate}

\begin{figure}[h]
\begin{center}
\includegraphics[width=\textwidth]{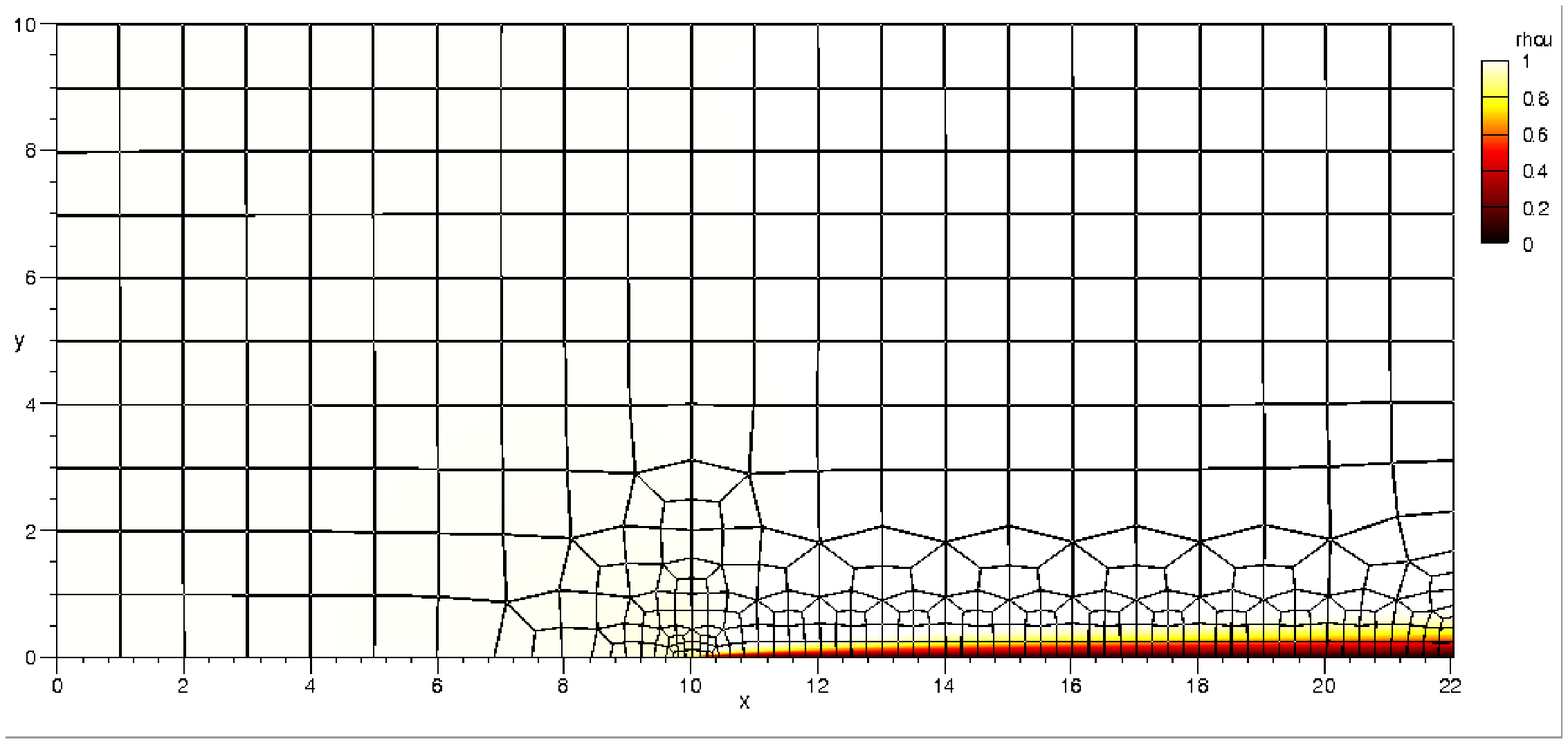}
\caption{Flat plate test case.}\label{fig:BLcontour}
\end{center}
\end{figure}

For this boundary layer test case, the mesh is constructed using 458 quadrilateral elements, and the simulations are computed with a Reynolds number of $\rm{Re}_{\infty}=6000$ (based on the reference length $L=12$) and a Mach number of M$_{\infty}=0.2$. Figure \ref{fig:BLcontour} shows the mesh and the distribution of the momentum in the x direction,$\rho u$.\\

On the boundary at $x=0$, a uniform inflow boundary condition was imposed. On the boundary $y=0$, $x<10$, a free-slip boundary condition was prescribed, whereas for $y=0$, $x \ge 10$, a no-slip adiabatic wall boundary condition emulates the effect of the flat plate. On the remaining boundaries we use a subsonic outflow boundary condition where only the far-field pressure is specified.\\

\subsubsection{Multigrid method}\label{sec:FlatPlateMG}
\begin{figure}[h]
\begin{center}

\subfigure[Residual norm vs. iterations.]{\label{fig:PerformanceFASiter} \includegraphics[width=0.45\textwidth]{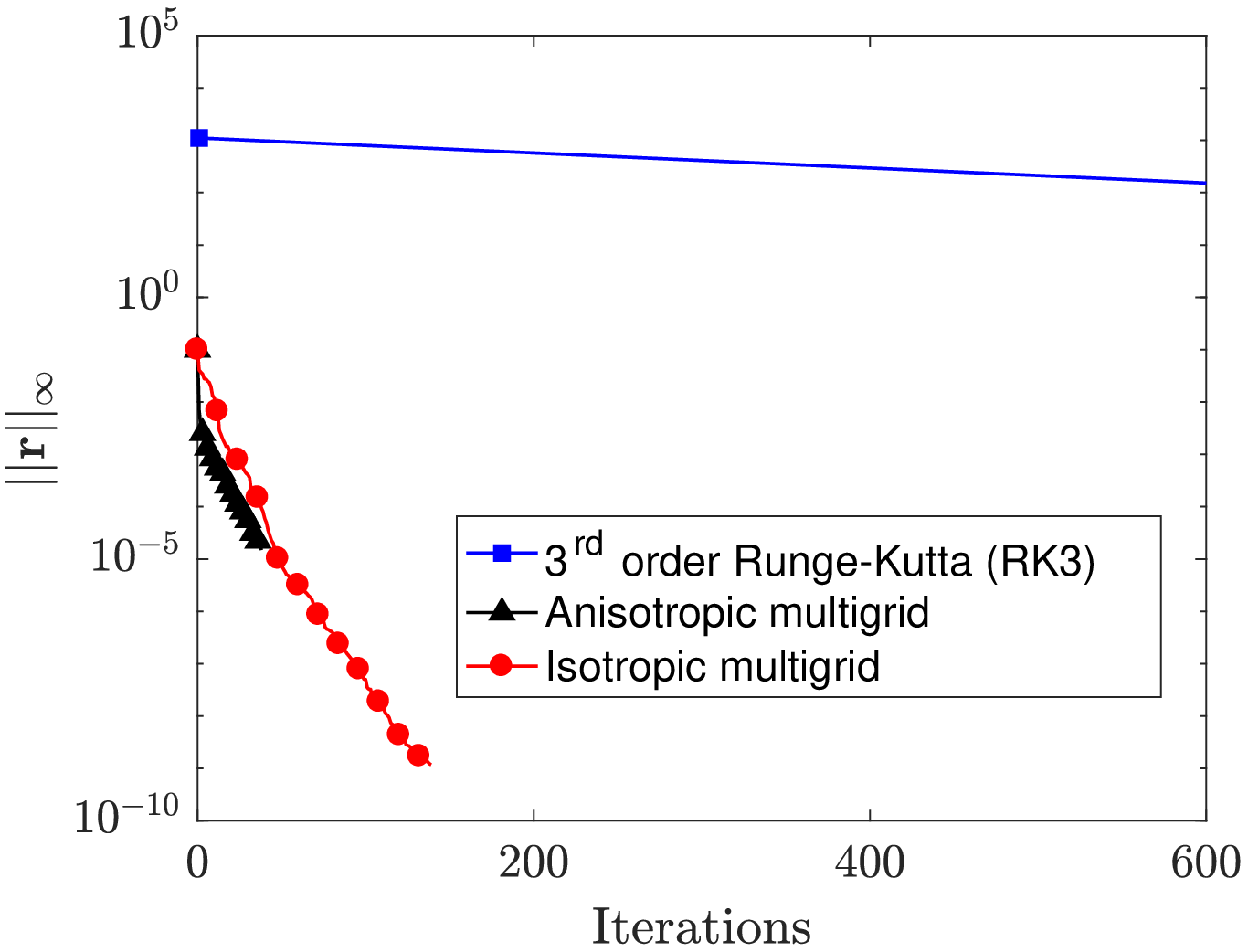}} \qquad
\subfigure[Residual norm vs. CPU-time.]{\label{fig:PerformanceFAS} \includegraphics[width=0.45\textwidth]{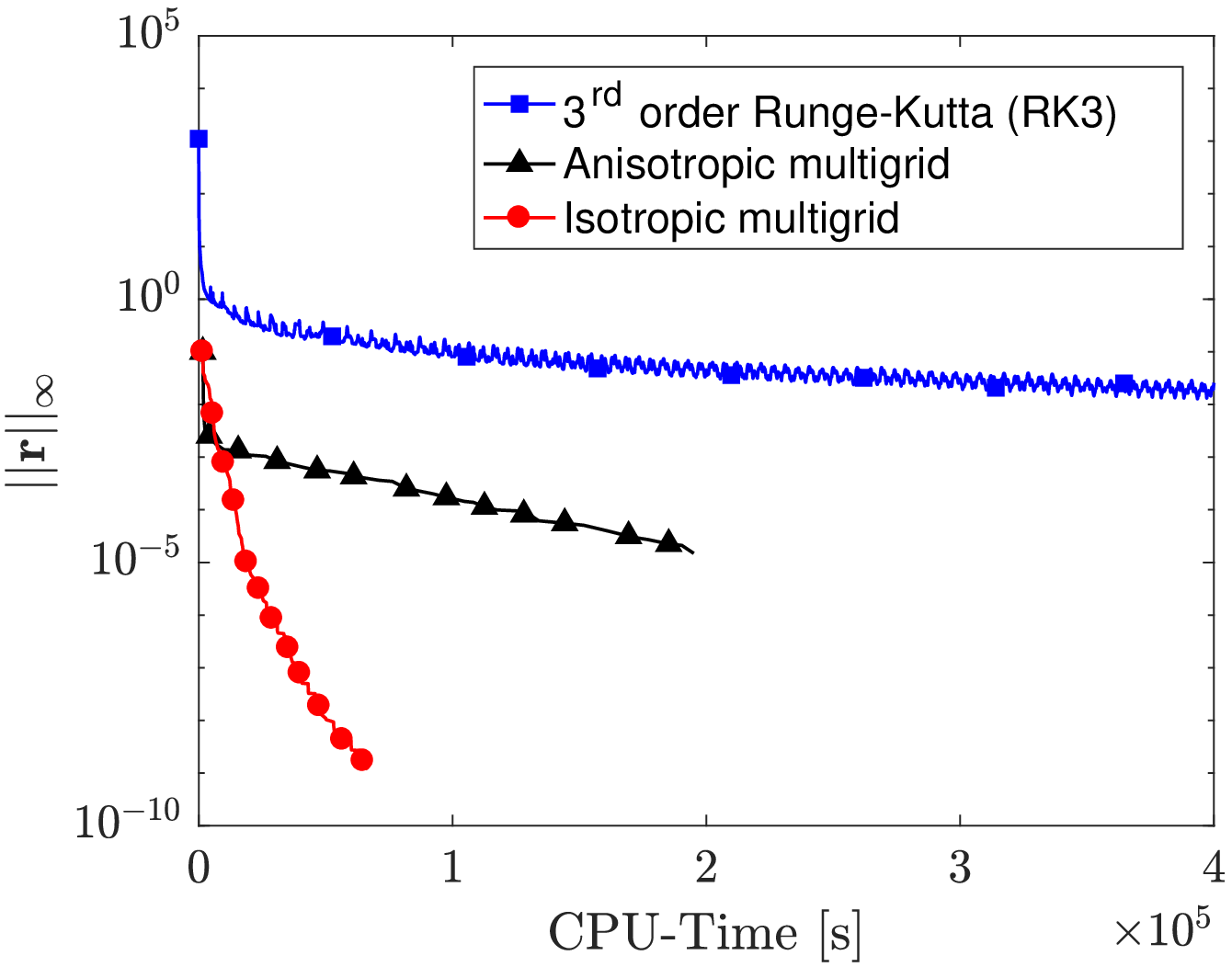}}
\caption{Comparison of performance of a classic RK3 method and a p-multigrid FAS (with a RK3 smoother) method for solving a subsonic boundary layer test case.}\label{fig:FASSpeedUp}
\end{center}
\end{figure}

First, a uniform mesh of order $N_1=N_2=10$ was simulated using different solution procedures. Figure \ref{fig:FASSpeedUp} shows the infinity norm of the residual (equation \eqref{eq:nonlinRes}) as a function of the iterations and the simulation time for the classic 3$^{rd}$ order Runge-Kutta scheme (RK3 - in blue), an isotropic p-multigrid FAS procedure (in red), and an anisotropic p-multigrid FAS procedure (in black), both of the latter using RK3 as smoother. All the results were obtained using $\beta_1^0=100$ pre-smoothing sweeps, $\beta_2^0=100$ post-smoothing sweeps, 400 smoothing sweeps on the coarsest multigrid level (common strategy for getting a good low-frequency representation), the smoothing tuning explained in section \ref{sec:TuninSmoothing}, and an FMG cycling strategy for obtaining an appropriate initial condition. As stated in section \ref{sec:MGcycling}, a residual-based strategy is used to control when the polynomial order is increased in the FMG cycle. A fixed residual of $10^{-1}$ must be obtained before the polynomial order is raised to the next FMG level. This threshold was selected because it showed good performance in preliminary tests.\\

It can be seen that the convergence rate of the multigrid strategies is much higher than the one of the completely explicit RK3 time integration, both in number of iterations and in CPU-time, even when the multigrid methods use the same RK3 as smoother. Additionally, the isotropic and anisotropic multigrid methods have a similar convergence rate with respect to the number of iterations, being the latter slightly better. However, when comparing the simulation times, it is remarkable that the isotropic FAS multigrid procedure is more efficient than the anisotropic one. The reason is that in the isotropic multigrid the lower multigrid levels have fewer degrees of freedom than in its anisotropic counterpart because the coarsening is done in all directions. For this reason, in next sections the anisotropic FAS will only be used as the $\hat \tau$-estimator (although during the estimation it is also used as a smoother) and the isotropic FAS will be mainly used as solver.

\subsubsection{Single-stage adaptation}
In this section, we study the computational cost involved in solving the boundary layer test case for different accuracy levels. To do that, we compare the results obtained using uniform adaptation with the ones obtained using the single-stage adaptation algorithm of section \ref{sec:NewPolorders}. In every case, the two main solvers considered in this paper were analyzed (the RK3 and the FAS solver). The single-stage adaptation process was performed for $N_{max}=10$ and $N_{max}=5$, where a reference mesh of $P_1=P_2=4$ was used. Notice that the use of such a coarse mesh as reference mesh is now possible because of the extrapolation capacities of the new estimation algorithm \cite{RuedaRamirez2018}. After adapting the mesh, the polynomial order jump across faces is limited to 

\begin{equation} \label{eq:poljump2D}
|N^{+}_i -N^{-}_i | \le 1,
\end{equation}
where the symbols $+$ and $-$ indicate the polynomial order in the direction $i$ of an element and its neighbor, respectively (the relative rotation between neighboring elements is taken into account). This condition provides robustness to the adapted mesh and is comparable with the \textit{two-to-one rule} that is usually employed in h-adaptation methods \cite{Liszka1997,Burgess2011}. Since the anisotropic truncation error estimator (equation \eqref{eq:AnisTruncError}) has been shown to generate more accurate extrapolations of the truncation error map than conventional $\hat \tau$-estimators \cite{RuedaRamirez2018}, the single-stage p-adaptation method (that is used in all cases) employs a 3V anisotropic V-cycle for estimating the \textit{isolated} truncation error, even when the time-marching solver is RK3.\\

A higher-order solution of order $N_1=N_2=15$ was used to estimate the relative error in the drag coefficient:

\begin{equation}
e_{\textit{drag}}^{N=15}= \frac{|C_d-C_d^{N=15}|}{C_d^{N=15}},
\end{equation}
where $C_d^{N=15}=0.211$, a value that is comparable to results in the literatute \cite{Lynn1995} for a flat plate at Re$_{\infty}=6000$. \\

\begin{figure}[htbp]
\begin{center}
\subfigure[Drag error vs. number of DOFs.]{\label{fig:DragErrorDOFs} \includegraphics[width=0.46\textwidth]{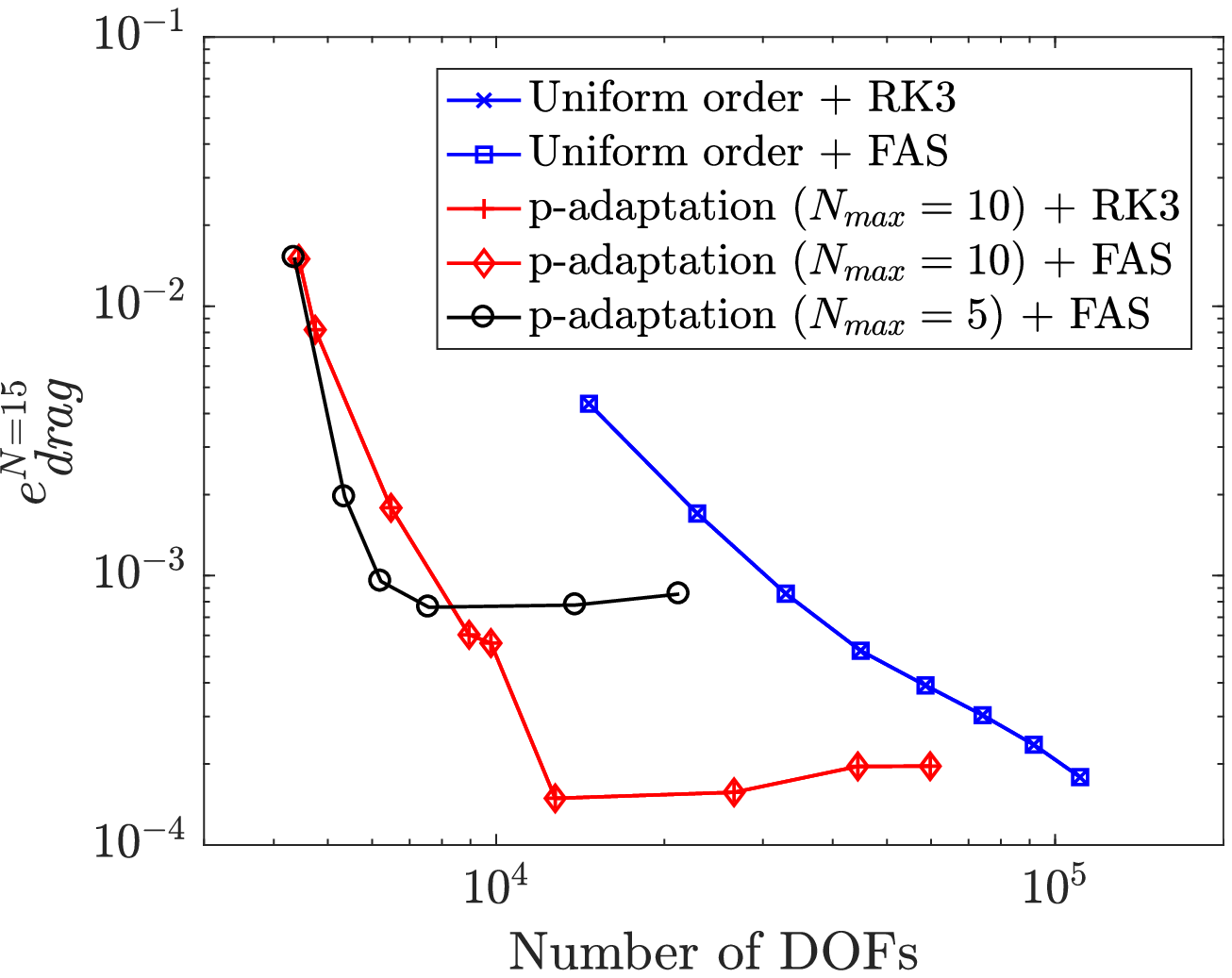}}\qquad
\subfigure[Drag error vs. CPU-Time.]{\label{fig:DragErrorTime} \includegraphics[width=0.46\textwidth]{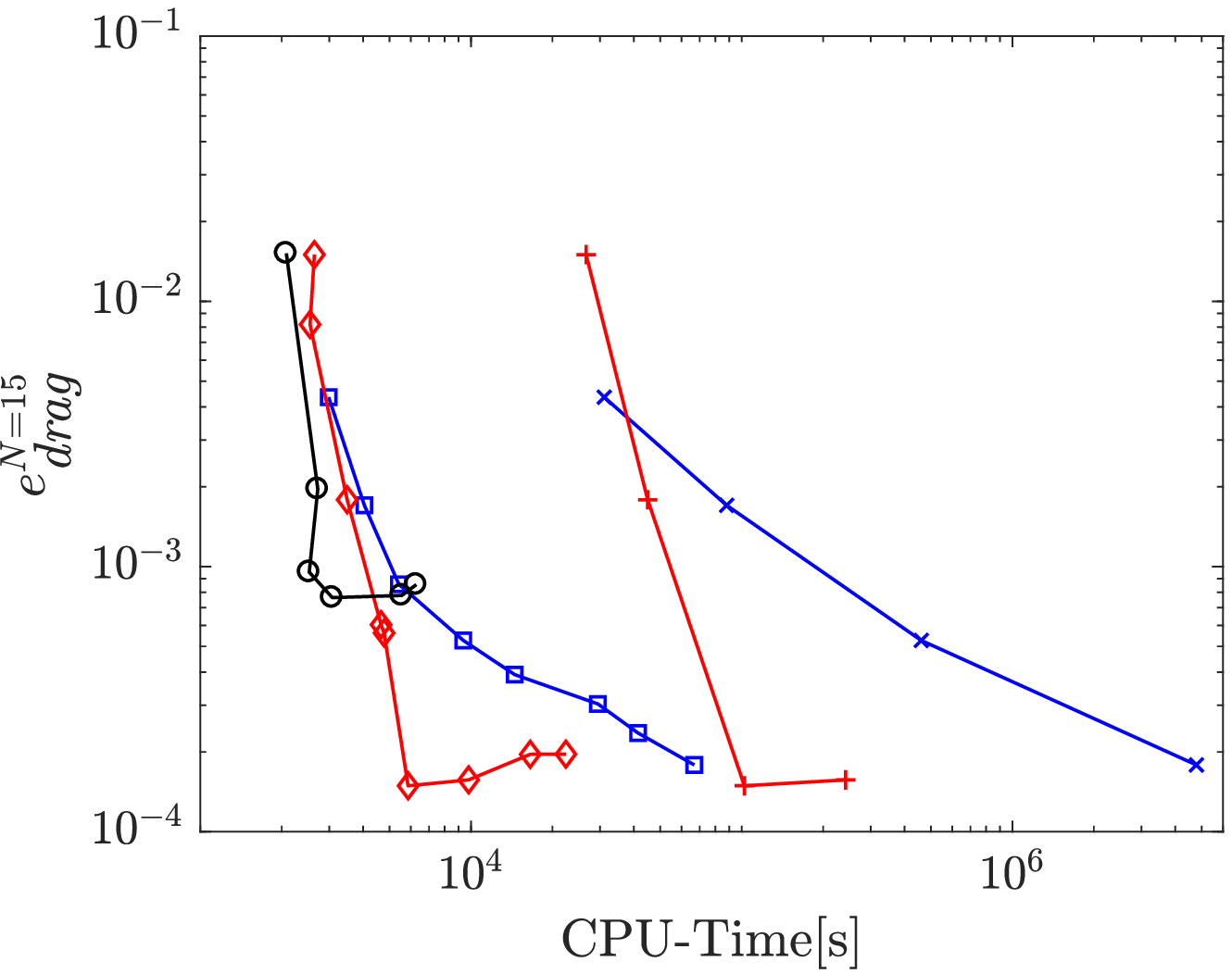}}

\caption{Relative error in the drag coefficient calculation for different methods. The reference drag $C_d^{P=15}$ was calculated on a uniformly refined mesh with $P = 15$. The blue lines represent uniform refinement, the red lines represent the $\hat \tau$-based p-adaptation procedure with $N_{max}=10$, and the black lines with $N_{max}=5$. Overlapping curves in (a).} \label{fig:DragError}
\end{center}
\end{figure}

\begin{figure}[h!]<
\begin{center}
\subfigure[Average polynomial order ($N_{av}$).]{\label{fig:AdaptedNav} \includegraphics[width=0.8\textwidth]{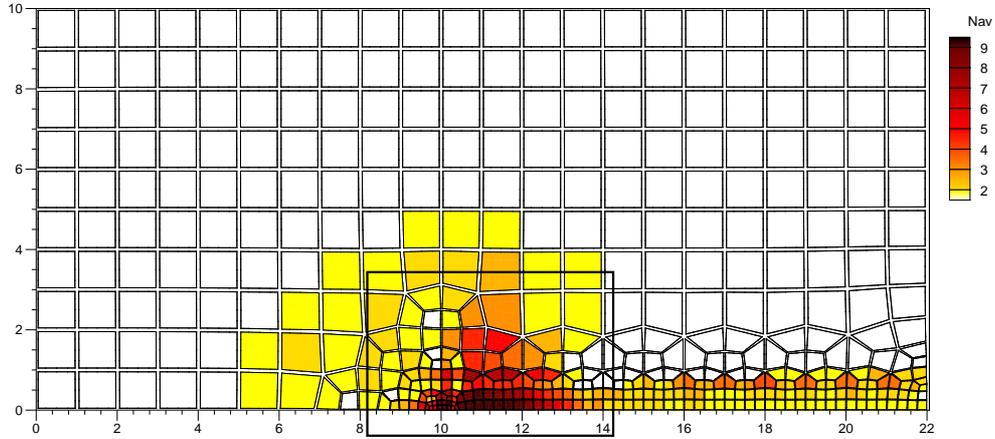}}\qquad
\subfigure[Detail of the Gauss-points.]{\label{fig:AdaptedNy} \includegraphics[width=0.7\textwidth]{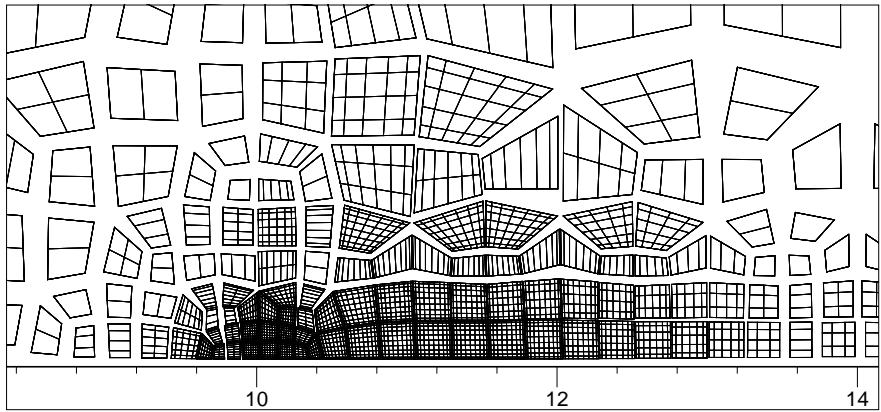}}

\caption{Contour indicating the final average polynomial orders after the adaptation procedure (a) and a detail of the Gauss-Points that shows the anisotropic nature of the p-adaptation method (b) for a threshold of $\hat \tau_{max}=10^{-3}$, which produces a relative drag error of $ e_{\textit{drag}}^{N=15} =1.49 \times 10^{-4}$. White boxes represent $N_1=N_2=1.$. $N_{av}=(N_1+N_2)/2$.
} \label{fig:AdaptedNxNy}
\end{center}
\end{figure}

The results obtained with the different methods are illustrated in Figure \ref{fig:DragError}. As shown in Figure \ref{fig:DragErrorDOFs}, the p-adapted meshes require a much fewer number of degrees of freedom for achieving a specific error than the uniformly adapted meshes. Note that the minimum relative error, that is achieved for low values of $\hat \tau_{max}$, tends to the relative error that corresponds to a mesh with uniform $N_{\textit{max}}$, in the same way as the minimum $\norm{\hat \tau}_{\infty}$ is a function of $N_{max}$ (see \cite{RuedaRamirez2018}). After meeting this plateau, no further improvement in the functional error is expected. This plateau is not necessarily obtained when all elements have $N_{max}$, as can be seen in Figure \ref{fig:AdaptedNxNy}. \\

As can be seen in Figure \ref{fig:DragErrorTime}, the p-adaptation procedures are especially efficient when high accuracy is needed. Using a low $N_{max}$ can lead to faster simulations, but the stagnation point is met sooner. It can also be observed that the most efficient procedure is the one that uses both FAS multigrid and p-adaptation. In fact, for the analyzed test case, this method achieves a better accuracy after a two hours of simulation than the classical approach (uniform order + RK3) after several days of computations. Table \ref{tab:SpeedUpOneStage} shows the CPU-time comparison of different solution procedures for reaching a drag error of at least $1.8 \times 10^{-4}$. The speed-up is as high as 815.76.\\

\begin{table}[h!]
\caption{Computation times and speed-up for the different methods for achieving a relative drag error of at least $1.8 \times 10^{-4}$ after converging until $\norm{\mathbf{r}}_{\infty} < 10^{-9}$.}
\begin{center}
\begin{tabular}{l|rrr}
Method  		  	& CPU-time[s] & Time [\%] & Speed-up \\ \hline \hline
RK3     			& $4.78\times 10^{6}$ & $100.00\%$ & $1.00$ \\ 
RK3 + p-adaptation 	& $1.02\times 10^{5}$ & $2.14\%$ & $46.72$ \\ 
FAS 				& $6.69\times 10^{4}$ & $1.40\%$ & $71.51$ \\ 
FAS + p-adaptation 	& $5.86\times 10^{3}$ & $0.12\%$ &$ 815.76$ \\ 
\end{tabular}
\label{tab:SpeedUpOneStage}
\end{center}
\end{table}

Figure \ref{fig:AdaptedNxNy} shows the final polynomial orders as computed by the proposed method for $\hat \tau_{max}=10^{-3}$ (equivalent to a drag error of $e_{\textit{drag}}^{N=15}=1.49 \times 10^{-4}$). It can be seen that intensive polynomial enrichment is performed on the leading edge of the flat plate around the singularity and on the regions where the boundary layer grows, as expected. Further polynomial enrichment can be observed in regions where the mesh size changes.\\

\subsubsection{Multi-stage adaptation} \label{sec:MultiStageRes2D}
In section \ref{sec:MultiStage}, we proposed a multi-stage p-adaptation procedure based on a full multigrid scheme with increasing polynomial orders and explored some of its theoretical advantages. Now, we apply this scheme to the boundary layer test case and analyze when it may be advantageous.\\

Moreover, in the last section we showed how the accuracy of the solution can be increased by choosing a higher $N_{max}$. Nonetheless, a higher $N_{max}$ represents a larger truncation error map. The calculations needed for generating the larger map are not computationally intensive \cite{RuedaRamirez2018}. However, as we increase the area of the map where we have to extrapolate the values, the uncertainty of the estimations also increases.\\

\begin{figure}[h]
\begin{center}
\includegraphics[width=0.5\textwidth]{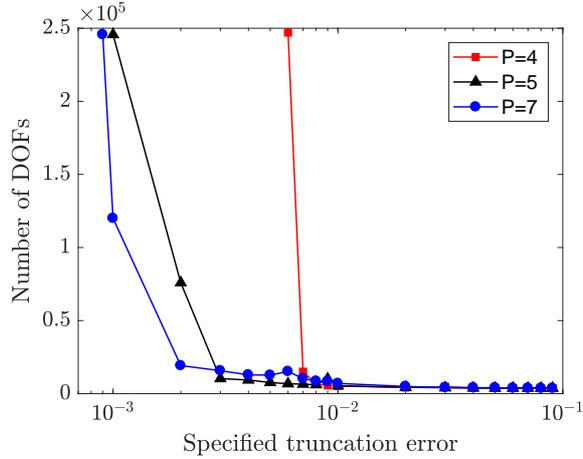}
\caption{Number of degrees of freedom obtained after adapting the mesh with different thresholds ($\hat \tau_{max}$) and different reference meshes ($P$) for $N_{max}=20$.}\label{fig:DOFsVsTauMax}
\end{center}
\end{figure}

Figure \ref{fig:DOFsVsTauMax} shows the number of degrees of freedom of the mesh after a single-stage adaptation procedure ($9\times 10^{-4} \le \hat \tau_{max} < 10^{-1}$ and $N_{max}=20$) for reference meshes of different order. As can be seen, the number of DOFs increases drastically when the specified truncation error is reduced below a certain value. This behavior occurs sooner for low-order reference meshes, where some elements are over-enriched to $N_{max}$. In fact, the polynomial order of the reference mesh is related to the maximum polynomial order it can accurately extrapolate the truncation error to. This relation is highly dependent on the PDE being approximated and the h-size of the mesh. Preliminary tests showed that, for the cases presented in this paper, a reference mesh of order $P$ can extrapolate accurately up to $2P$.\\

For high values of $N_{max}$, a multi-stage p-adaptation procedure becomes very useful. As was explained in section \ref{sec:MultiStage}, instead of starting with a high-order reference mesh (which can be very expensive), a coarse reference mesh of order $P=\mathcal{P}_1$ is chosen to estimate the truncation error. With the estimation, the regions where a low-order approximation is enough are identified. Afterwards, the p-adaptation algorithm sets the polynomial orders of the mesh according to the $\hat \tau$-estimation and limits the over-enrichment in more complex flow regions to $\mathcal{P}_2$. In the second adaptation process at $P=\mathcal{P}_2$, and in subsequent adaptation stages, the polynomial orders of the mesh are corrected with a more accurate error estimation at hand. \\

\begin{figure}[htbp]
\begin{center}
\subfigure[$N_{max}=20$.]{\label{fig:PerfNmax20} \includegraphics[width=0.46\textwidth]{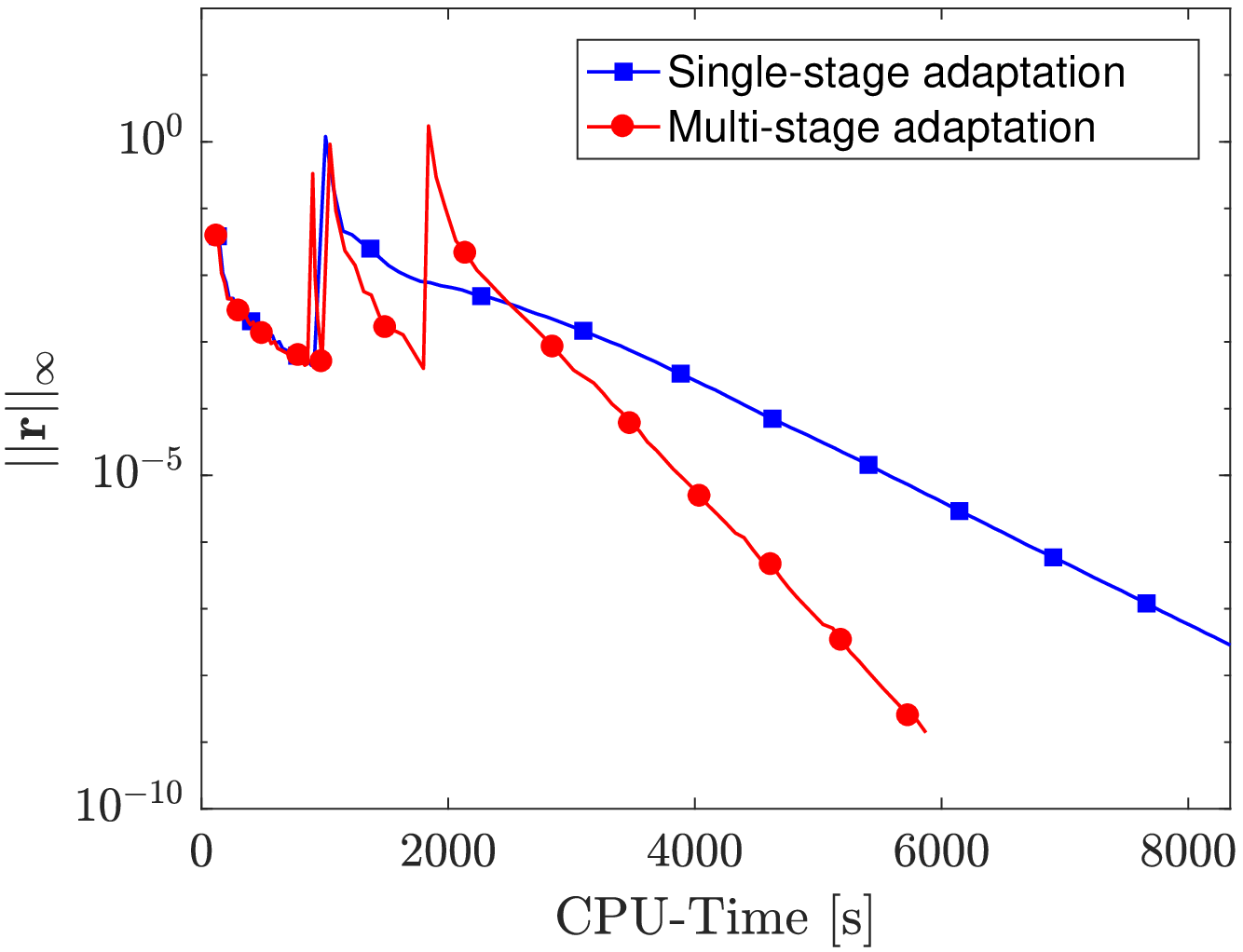}}\qquad
\subfigure[$N_{max}=30$.]{\label{fig:PerfNmax30} \includegraphics[width=0.46\textwidth]{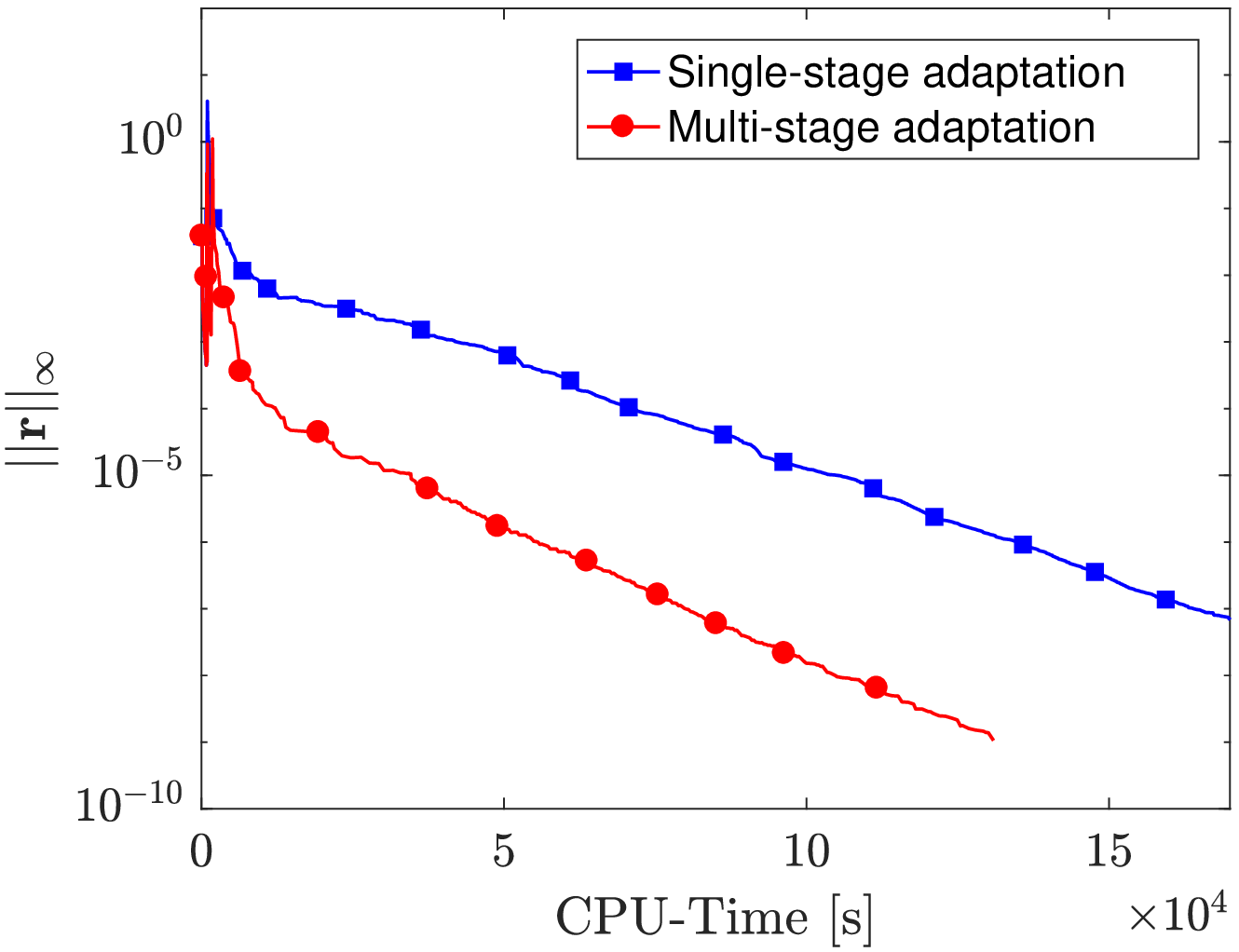}}

\caption{Comparison of a single-stage and a multi-stage adaptation process for solving the boundary layer test case with a truncation error threshold of $\hat \tau_{max}=4 \times 10 ^{-3}$: $N_{max}=20$ (a), and $N_{max}=30$ (b).} \label{fig:OneStageVsMultiStage}
\end{center}
\end{figure}

In order to illustrate how this method can reduce the computational cost of highly accurate simulations, we present a comparison of the convergence of the single-stage and the multi-stage p-adaptation procedures for $\hat \tau_{max}=4 \times 10^{-3}$, and $N_{max}=20$ (Figure \ref{fig:PerfNmax20}) and $N_{max}=30$ (Figure \ref{fig:PerfNmax30}). The reference meshes of the multi-stage algorithm were selected at $\mathcal{P}_1=4$, $\mathcal{P}_2=8$ and $\mathcal{P}_3=16$. The measured speed-up is 1.69 for $N_{max}=20$ and 1.72 for $N_{max}=30$ with respect to the single-stage adaptation.

\subsection{3D Flow around a Sphere} \label{sec:Sphere}

\begin{figure}[tbp]
\begin{center}
\includegraphics[width=.8\textwidth]{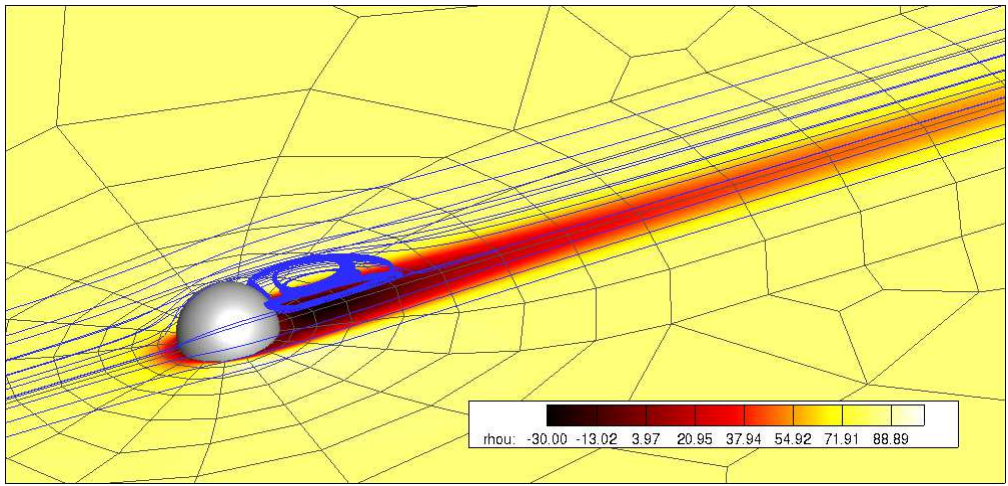}
\caption{Sphere at Re$_{\infty}=200$.}\label{fig:Spherecontour}
\end{center}
\end{figure}

For this test case, the mesh is constructed with 1904 hexahedral elements, and the simulations are computed with a Reynolds number of Re$_{\infty}=200$ and a Mach number of M$_{\infty}=0.2$. The curvilinear hexahedral mesh has a mapping order $M=3$ and was created using the HOPR package \cite{hindenlang2015mesh}. Figure \ref{fig:Spherecontour} shows the mesh and the distribution of the conserved variable $\rho u$ around the sphere.\\

In order to assess the properties of the representations obtained after performing $\hat \tau$-based adaptation, we use a relative drag error that is computed against a high-order solution of order $N=12$ in the same mesh:

\begin{equation}
e_{\textit{drag}}^{N=12}= \frac{|C_d-C_d^{N=12}|}{C_d^{N=12}}.
\end{equation}

Table \ref{tab:DragSphere} shows a comparison between the reference drag coefficient obtained in this work and in other studies. 

\begin{table}[h!]
\caption{Drag coefficient for sphere at Re$_{\infty}=200$.}
\begin{center}
\begin{tabular}{lr}
Author									&	Value 	\\ \hline\hline
Campregher et al. \cite{Campregher2009}	& 0.815 	\\
Fornberg \cite{Fornberg1988}			& 0.7683	\\
Fadlun et al. \cite{Fadlun2000}			& 0.7567	\\
This work								& 0.7771  \\
\end{tabular}
\label{tab:DragSphere}
\end{center}
\end{table}

\subsubsection{3D Considerations}

Since in the general 3D p-nonconforming DGSEM the mapping order in every direction is limited by the solution order as $M_i \le N_i/2$ (as indicated in section \ref{sec:DGSEM}), and considering that the p-adapted meshes are in general p-nonconforming, the minimum polynomial order after p-adaptation is set to $N_{min}=2$. Additionally, taking into account the observations made in section \ref{sec:HOcoarsening}, we use high-order coarsening and $N_{\textit{coarse}}=1$ for the p-multigrid method before and after p-adaptation. \\

Moreover, in order to represent the curved boundary on the sphere as exactly as possible, after the p-adaptaion, a conforming algorithm changes the polynomial orders of all elements on that surface, so that there is no polynomial order jump across their faces. This allows using a mapping of order $M_i \le N_i$ there. \\

Finally, let us remark that in 3D, the condition of equation \eqref{eq:poljump2D} (polynomial jump across faces of 1) can cause a steep increase in the number of degrees of freedom because the polynomial enriching is transmitted in three directions. Therefore, for this test case the polynomial order jump across faces after p-adaptation is softened to 

\begin{equation}
N^{+}_i \ge \bigg\lfloor \frac{2}{3} N^{-}_i \bigg\rfloor,
\end{equation}
where $\lfloor \cdot \rfloor$ is the integer part floor function.\\

This condition showed to provide enough robustness to the p-adapted representations and lowered the number of degrees of freedom of the adapted meshes. The conforming algorithm that is used on the sphere boundary and the algorithm that controls the polynomial order jump everywhere must be executed iteratively, until no further changes are needed, to ensure that the final mesh has all the desired properties.

\subsubsection{Single-stage adaptation}
The single-stage adaptation process is performed for $N_{max}=7$, where a reference mesh of order $P_1=P_2=P_3=5$ is used. Different values of the specified truncation error threshold were tested in the range $10^{-1} \le \hat \tau_{max} \le 10^{-4}$.\\

The isotropic and conforming reference mesh is iterated down to a residual of $\norm{\mathbf{r}}_{\infty} \le \hat \tau_{max}/10$ using a p-multigrid algorithm with $\beta_1^0=\beta_2^0=100$ pre- and post-smoothing sweeps, and 400 smoothing sweeps on the coarsest multigrid level. After the p-adaptation, the pre- and post-smoothing sweeps are $\beta_1^0=\beta_2^0=50$, and the  number of smoothing sweeps on the coarsest multigrid level is 200. This combination exhibited the best performance. The smoothing tuning detailed in section \ref{sec:TuninSmoothing}  is used and an FMG cycling strategy is employed for obtaining an appropriate initial condition with a residual of $\norm{\mathbf{r}}_{\infty} \le 1.0$. \\

\begin{figure}[htbp]
\begin{center}
\subfigure[Drag error vs. number of DOFs.]{\label{fig:SphereDragDOFs} \includegraphics[width=0.45\textwidth]{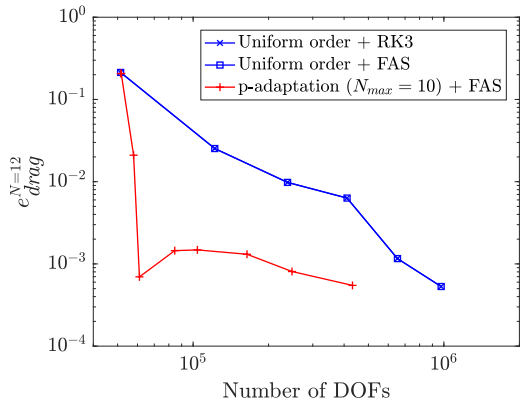}}\qquad
\subfigure[Absolute lift vs. number of DOFs.]{\label{fig:SphereLiftDOFs} \includegraphics[width=0.45\textwidth]{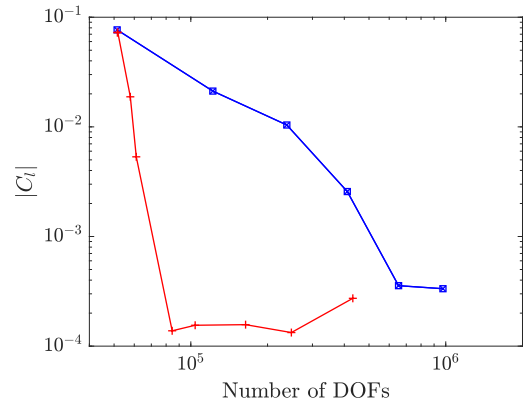}}

\subfigure[Drag error vs. CPU-Time.]{\label{fig:SphereDragTime} \includegraphics[width=0.45\textwidth]{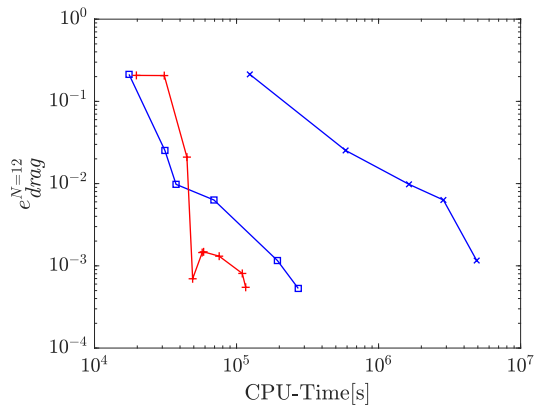}}
\subfigure[Absolute lift vs. CPU-Time.]{\label{fig:SphereLiftTime} \includegraphics[width=0.45\textwidth]{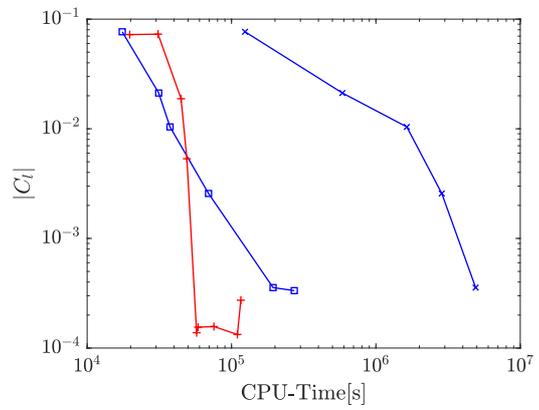}}

\caption{Relative error in the drag and lift coefficients for different methods on the sphere. The blue lines represent uniform refinement, and the red lines represent the $\hat \tau$-based p-adaptation procedure with $N_{max}=7$.} \label{fig:SphereErrors}
\end{center}
\end{figure}

The relative drag error and the absolute lift of the adapted meshes are assessed. Figure \ref{fig:SphereErrors} shows a comparison between the errors obtained using the $\hat \tau$-based adaptation procedure and the ones using uniform p-refinement. As can be observed in Figures \ref{fig:SphereDragDOFs} and \ref{fig:SphereLiftDOFs}, the number of degrees of freedom is greatly reduced for the same accuracy when using the $\hat \tau$-based p-adaptation. The maximum error in both coefficients is related to the one obtained with a uniform mesh of $N_1=N_2=N_3=N_{min}=2$, as expected. Similarly, the minimum error tends to stagnate at a value that is comparable to the one obtained for $N_1=N_2=N_3=N_{min}=7$. \\

It is interesting to notice that, for high truncation error thresholds, the $\hat \tau$-based adaptation does not provide an advantage in CPU-time (Figures \ref{fig:SphereDragTime} and \ref{fig:SphereLiftTime}). This is due to the cost of obtaining a semi-converged solution on the reference mesh of $P=5$, for cases where the final polynomial order post-adaptation is $N<5$. Additionally, let us remark that the rate of convergence in CPU-time deteriorates after the p-adaptation. This is because the p-anisotropic nonconforming representations are \textit{more difficult} to solve. Further investigation on multigrid, or other solution methods, could improve the speed-ups here observed.

Using the data provided by Figures \ref{fig:SphereDragTime} and \ref{fig:SphereLiftTime}, it is possible to compute the speed-up as a function of the drag or lift errors. Table \ref{tab:SpeedUpSphere} shows the computation times and speed-ups achieved for the lowest error obtained in lift and drag ($e_{\textit{drag}}^{N=12}$ and $|C_l|$). The maximum speed-up is 151.94 for this 3D challenging case.\\ 

\begin{figure}[h!]
\begin{center}
\subfigure[Average polynomial order ($N_{av}$).]{\label{fig:SphereNav} \includegraphics[width=0.7\textwidth]{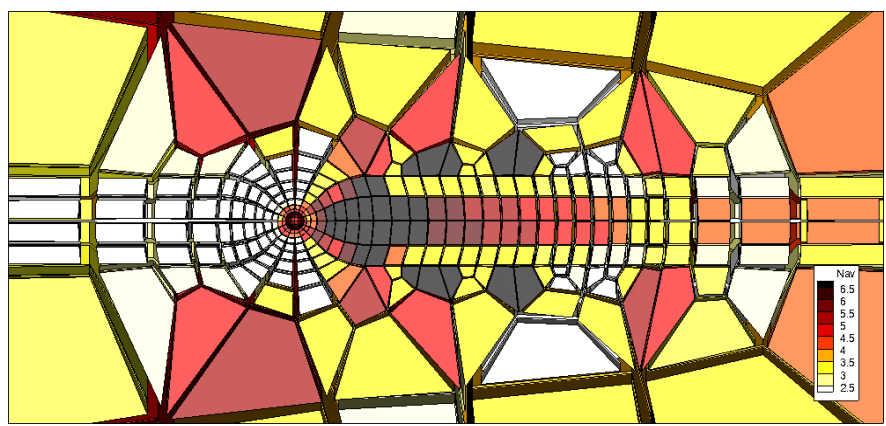}}\qquad
\subfigure[Detail of the Gauss-points.]{\label{fig:SphereDetail} \includegraphics[width=0.7\textwidth]{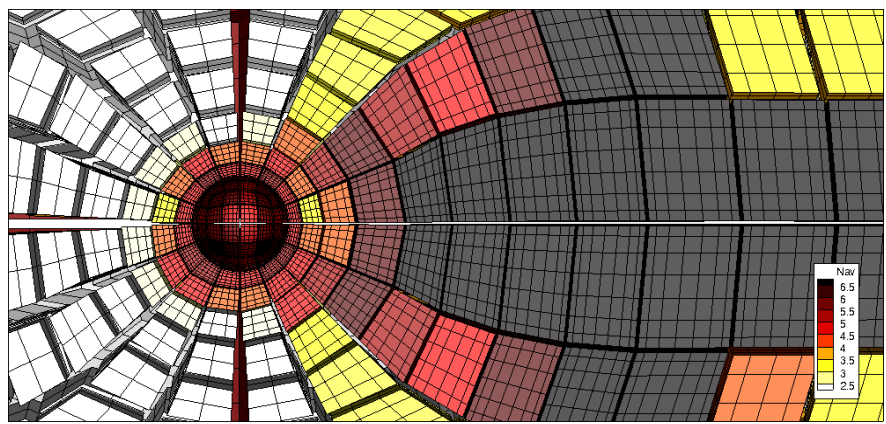}}


\caption{Contours indicating the final polynomial orders after p-adaptation for the sphere test case: Average polynomial orders ($N_{av}$) (a) and a detail of the Gauss-Points that shows the anisotropic nature of the p-adaptation method (b) for a threshold of $\hat \tau_{max}=4 \times 10^{-4}$. White boxes represent $N_1=N_2=N_3=2.$. $N_{av}=(N_1+N_2+N_3)/3$.
} \label{fig:AdaptedSphere}

\end{center}
\end{figure}

Figure \ref{fig:AdaptedSphere} illustrates the polynomial order distribution after p-adaptation for $\hat \tau_{max}=4\times 10^{-4}$, which corresponds to a drag error of $e_{\textit{drag}}=8.08\times 10^{-4}$ and an absolute lift of $|C_l|=1.33\times 10^{-4}$. It can be seen that intensive polynomial enrichment is performed on the recirculation bubble, the wake, and on the boundary layer, as expected. Further polynomial enrichment can be observed in regions where the mesh size changes drastically. In particular, we observe that the polynomial enrichment is higher on the wake that on the recirculation bubble because of the element sizes of the available mesh. \\

\begin{table}[h!]
\caption{Computation times and speed-up for the different methods after converging until $\norm{\mathbf{r}}_{\infty} < 10^{-9}$}
\begin{center}
\begin{tabular}{l|rrr|rrr|}
		& \multicolumn{3}{c|}{Drag coefficient ($e_{\textit{drag}} \le \times 5.31 \times 10^{-4}$)} & \multicolumn{3}{c|}{Lift coefficient ($|C_l| \le 3.34 \times 10^{-4}$)} \\
Method  & CPU-time[s] & Time [\%] & Speed-up & CPU-time[s] & Time [\%] & Speed-up \\ \hline \hline
RK3     & $7.46\times 10^{6}$ & $100.00\%$ & $1.00$ 
		& $7.46\times 10^{6}$ & $100.00\%$ & $1.00$  \\ 
FAS 	& $2.72\times 10^{5}$ & $3.65\%$   & $27.41$  
		& $2.72\times 10^{5}$ & $3.65\%$   & $27.41$ \\ 
FAS + p-adaptation & $4.91\times 10^{4}$ & $0.68\%$   & $151.94$ 
				   & $5.80\times 10^{4}$ & $1.06\%$   & $128.55$ \\ 
\end{tabular}
\label{tab:SpeedUpSphere}
\end{center}
\end{table}


\subsubsection{Multi-stage adaptation}
\begin{figure}[h]
\begin{center}
\includegraphics[width=0.5\textwidth]{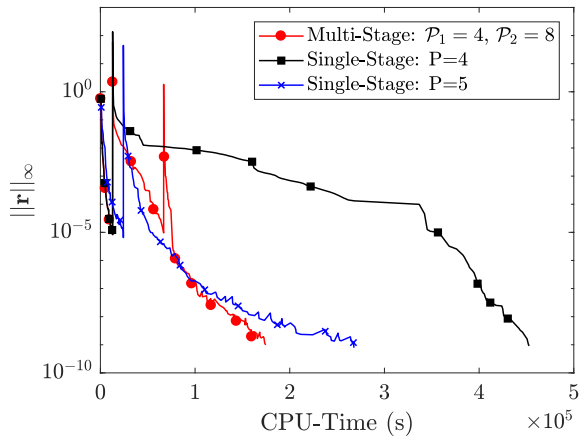}
\caption{Comparison of single-stage ($P=4$ and $P=5$) and multi-stage adaptation ($\mathcal{P}_1=4$, $\mathcal{P}_2=8$) processes for the sphere. $N_{\textit{max}}=11$, $\hat \tau_{\textit{max}}=10^{-4}$.}\label{fig:SphereMultiStage}
\end{center}
\end{figure}

The multi-stage adaptation procedure introduced in section \ref{sec:MultiStage} becomes useful when the maximum allowable polynomial order after adaptation ($N_{\textit{max}}$) is increased and the specified \textit{isolated} truncation threshold ($\hat \tau_{\textit{max}}$) is low. In this section, we use the multi-stage p-adaptation procedure on the sphere test case and set the maximum polynomial order after adaptation to $N_{\textit{max}}=11$, the truncation error threshold to $\hat \tau_{\textit{max}}=10^{-4}$, and the adaptation stages at $\mathcal{P}_1=4$ and $\mathcal{P}_2=8$. Figure \ref{fig:SphereMultiStage} shows a comparison of performance (in CPU-Time) between the multi-stage p-adaptation procedure and two single-stage procedures with $P=4$ and $P=5$. The maximum polynomial order after adaptation in the single stage cases is also $N_{\textit{max}}=11$.\\

As can be observed, the convergence rate (with respect to CPU-time) of the multi-stage p-adapted mesh is higher than for the single-stage p-adapted meshes. The reason is that the former has fewer degrees of freedom. Table \ref{tab:MultiStageSphere} shows a summary of results for the simulations of this section.

\begin{table}[h!]
\caption{Summary of performance for single- and multi-stage simulations with $\hat \tau_{\textit{max}}=10^{-4}$}
\begin{center}
\begin{tabular}{l|rrrrrr}
Adaptation strategy  &	DOFs(1) & DOFs(2) & CPU-Time(s) & Speed-up & $e_{\textit{drag}}$ & $|C_l|$ \\ \hline\hline
Single-Stage: $P=4$  & $1.07 \times 10^6$ & --- & $4.53 \times 10^5$ & 1.00 & $2.27 \times 10^{-5}$ 	& $1.39 \times 10^{-5}$ \\
Single-Stage: $P=5$	 & $7.90 \times 10^5$ & --- & $2.68 \times 10^5$ & 	     1.69 & $ 3.57 \times 10^{-5}$ 	& $1.90 \times 10^{-4}$\\
Multi-Stage: $\mathcal{P}_1=4$, $\mathcal{P}_2=8$ 
					 & $6.20 \times 10^5$ & $3.85 \times 10^5$
					 & $1.75 \times 10^5$ & 2.59 & $4.50 \times 10^{-5}$
					 & $2.12 \times 10^{-5}$\\
\end{tabular}
\label{tab:MultiStageSphere}
\end{center}
\end{table}

The number of degrees of freedom for the single-stage $P=4$ is the highest, since in that case many elements are enriched to the maximum $N_1=N_2=N_3=11$ due to problems in the error estimation (as explained in section \ref{sec:MultiStageRes2D}). In the single-stage $P=5$ this behavior is also observed, but to a lesser extent. On the other hand, in the multi-stage case the number of degrees of freedom in the first stage is limited by the condition $N_{\textit{max},1}=8$, and the distribution of polynomial orders is then corrected in the second stage, where the number of degrees of freedom decreases, even though the maximum polynomial order is $N_{\textit{max},2}=11$. It is remarkable that the multi-stage adapted mesh can achieve comparable drag and lift errors with about one third of the degrees of freedom and a speed-up of 2.59 with respect to the single-stage $P=4$.

\section{Conclusions} \label{sec:Conclusions}
In this paper, we have developed a coupled solver using truncation error estimators, anisotropic p-adaptation and multigrid. The most important conclusions of this work are:
\begin{enumerate}
\item A novel anisotropic p-adaptation multigrid algorithm is presented which uses the multigrid method both as a solver and as a truncation error estimator.

\item The coupling of single-stage p-adaptation strategies and multigrid methods resulted in a speed-up of 816 for a 2D boundary layer case and of 152 for the 3D sphere case.

\item The technique for evaluating the truncation error by Rueda-Ramírez et al. \cite{RuedaRamirez2018} can be performed directly inside an anisotropic multigrid procedure needing only a few additional operations.

\item Isotropic multigrid methods show better performance than anisotropic multigrid methods. The reason is that the successive coarse grids are cheaper to compute when the polynomial order is reduced in all directions.

\item A multi-stage p-adaptation technique based on coupling $\tau$-estimations and multigrid was developed. Experiments show that multi-stage is advantageous for highly accurate simulations compared with single-stage adaptation procedures. The multi-stage procedure showed to be a promising alternative for 3D simulations, since coarser reference meshes can be used: the elements that do not need to be enriched are identified early and their polynomials are frozen in a low value. The achieved speed-ups with this methods were as high as 2.59 with respect to the single-stage adaptation.

\end{enumerate}

\section*{Acknowledgments}
The authors would like to thank David Kopriva for his friendly advise and cooperation. This project has received funding from the European Union’s Horizon 2020 Research and Innovation Program under the Marie Skłodowska-Curie grant agreement No 675008 for the SSeMID project.\\

The authors acknowledge the computer resources and technical assistance provided by the \textit{Centro de Supercomputación y Visualización de Madrid} (CeSViMa).

\section*{References}
\bibliography{Biblio.bib}

\appendix

\section{The Navier-Stokes equations}\label{sec:NS}

The compressible Navier-Stokes equations in conservative form can be written in non-dimensional form as
\begin{equation}
\mathbf{q}_t + \nabla \cdot \left( \mathscr{F}^a - \mathscr{F}^{\nu} \right) = \mathbf{s}, 
\end{equation}

where the conserved variables are $\mathbf{q} = \left( \rho, \rho u, \rho v, \rho w, \rho e  \right)^T$, $\mathbf{s}$ is an external source term, and $\mathscr{F}^a$ and $\mathscr{F}^{\nu}$ are called the advective and diffusive flux tensors, respectively, which depend on $\mathbf{q}$. Expanding the fluxes in Cartesian coordinates leads to the expression,

\begin{equation}
\mathbf{q}_t + \mathbf{f}^a_x + \mathbf{g}^a_y + \mathbf{h}^a_z - \frac{1}{\rm{Re}_{\infty}} \left( \mathbf{f}^{\nu}_x + \mathbf{g}^{\nu}_y + \mathbf{h}^{\nu}_z \right) = \mathbf{s}.
\end{equation}
Here, $\textrm{Re}_{\infty}=V_{\infty} L_{\infty} \rho_{\infty}/ \mu_{\infty}$ is the Reynolds number in the far-field. The advective fluxes are then defined as

\begin{equation}
\mathbf{f}^a =
          \begin{bmatrix}
           \rho u      \\           
           p + \rho u^2\\
           \rho u v    \\
           \rho u w    \\
           u (\rho e + p)
          \end{bmatrix},
          		\mathbf{g}^a = 
        				\begin{bmatrix}
				           \rho v      \\           
				           \rho u v    \\
           				   p + \rho v^2\\
				           \rho v w    \\
				           v (\rho e + p)
        				\end{bmatrix},
          						\mathbf{h}^a = 
          								\begin{bmatrix}
								           \rho w      \\           
								           \rho u w    \\
				           				   \rho v w    \\
								           p + \rho w^2\\
								           w (\rho e + p)
			          					\end{bmatrix},
\end{equation}
where the pressure $p$ is computed using the calorically perfect gas approximation. On the other hand, the diffusive fluxes are defined as
\begin{align}
\mathbf{f}^{\nu} &=
          \begin{bmatrix}
           0           \\           
           \tau_{xx}   \\
           \tau_{xy}   \\
           \tau_{xz}   \\
           u \tau_{xx} + v \tau_{xy} + w \tau_{xz} + \frac{\kappa}{(\gamma - 1) \rm{Pr}_{\infty} \rm{M}_{\infty}^2} T_x 
          \end{bmatrix}, \\
          		\mathbf{g}^{\nu} &= 
        				\begin{bmatrix}
				           0           \\           
				           \tau_{yx}   \\
				           \tau_{yy}   \\
				           \tau_{yz}   \\
				           u \tau_{yx} + v \tau_{yy} + w \tau_{yz} + \frac{\kappa}{(\gamma - 1) \rm{Pr}_{\infty} \rm{M}_{\infty}^2} T_y \\
        				\end{bmatrix}, \\
          						\mathbf{h}^{\nu} &= 
          								\begin{bmatrix}
								           0           \\           
								           \tau_{zx}   \\
								           \tau_{zy}   \\
								           \tau_{zz}   \\
								           u \tau_{zx} + v \tau_{zy} + w \tau_{zz} + \frac{\kappa}{(\gamma - 1) \rm{Pr}_{\infty} \rm{M}_{\infty}^2} T_z
			          					\end{bmatrix},
\end{align}
where $T$ is the temperature, $T_i$ its spatial derivatives, $\gamma$ is the heat capacity ratio, and $\kappa$ is the thermal diffusivity. The nondimensional far-field parameters are the Prandtl number, $\rm{Pr}_{\infty}=c_p \mu_{\infty} / \kappa_{\infty}$; and the Mach number, $\rm{M}_{\infty}=\norm{\mathbf{v}}/c$. The stress tensor components are computed using the Stokes hypothesis,
\begin{align}
\tau_{ij} &= \mu \left( \frac{\partial v_i}{\partial x_j} + \frac{\partial v_j}{\partial x_i} \right), i \neq j \\
\tau_{ii} &= 2 \mu \frac{\partial v_i}{\partial x_i} + \lambda \nabla \cdot \mathbf{v},
\end{align}
with $\mu$ the fluid's viscosity, $\lambda=-\frac{2}{3} \mu$ the bulk viscosity coefficient, and $\mathbf{v}$ the flow velocity. For the simulations in this paper we chose the typical parameters for air: $\rm{Pr} = 0.72$, $\gamma = 1.4$, while $\mu$ and $\kappa$ are calculated using Sutherland's law.

\end{document}